\documentclass[aos]{imsart}

\RequirePackage{amsthm,amsmath,amsfonts,amssymb}
\RequirePackage[authoryear]{natbib}
\RequirePackage[colorlinks,citecolor=blue,urlcolor=blue]{hyperref}
\RequirePackage{graphicx}

\startlocaldefs
\theoremstyle{plain}

\newtheorem{theorem}{Theorem}[section]
\theoremstyle{remark}
\newtheorem{definition}[theorem]{Definition}
\newtheorem*{example}{Example}

\usepackage{amsmath,amssymb,amsfonts,amsthm,bbm,empheq}
\usepackage{epic,eepic,epsfig,longtable}
\usepackage{multirow,multicol,verbatim}
\usepackage{array}
\usepackage{lscape}
\usepackage{graphicx}
\usepackage{latexsym}
\usepackage{comment}
\usepackage{booktabs}
\usepackage{epsfig}
\usepackage{CJK}
\usepackage{color}
\usepackage{mathtools}
\usepackage{extarrows}
\usepackage{centernot}
\usepackage{caption}
\usepackage{subcaption}
\usepackage{tikz}
\usetikzlibrary{shapes.multipart,shapes.geometric,spy}

\numberwithin{equation}{section}

\NewCommandCopy{\baccent}{\b}
\AtBeginDocument{%
  \DeclareRobustCommand{\b}[1]{\ifmmode\mathbf{#1}\else\baccent{#1}\fi}%
}

\renewcommand{\emptyset}{\varnothing}

\def\[{\left [}  \def\]{\right ]} \def\({\left (}  \def\){\right )}

\makeatletter
\def\underbar#1{\underline{\sbox\tw@{$#1$}\dp\tw@\z@\box\tw@}}
\makeatother

\def\tilde{\widetilde}

\def\newpage{\vfill\eject}
\def\today{\ifcase\month\or
  January\or February\or March\or April\or May\or June\or
  July\or August\or September\or October\or November\or December\fi
  \space\number\day, \number\year}

\usepackage{stackengine} 
\stackMath

\setstackgap{L}{8pt}

\allowdisplaybreaks
\def\thesection{\arabic{section}}
\usepackage{verbatim}

\usepackage{lipsum}
\usepackage{enumitem}
\usepackage{hyperref}
\usepackage{xcolor}

\def\:={\coloneqq}
\def\=:{\eqqcolon}
\theoremstyle{plain}
\ifx\theorem\undefined
  \newtheorem{theorem}{Theorem}[section]
\fi
\ifx\proposition\undefined
  \newtheorem{proposition}{Proposition}[section]
\fi
\ifx\corollary\undefined
  \newtheorem{corollary}{Corollary}[section]
\fi
\ifx\lemma\undefined
  \newtheorem{lemma}{Lemma}[section]
\fi

\theoremstyle{definition}
\ifx\definition\undefined
  \newtheorem{definition}{Definition}[section]
\fi
\ifx\assumption\undefined
  \newtheorem{assumption}{Assumption}
\fi
\ifx\example\undefined
  \newtheorem{example}{Example}[section]
\fi
\ifx\exercise\undefined
  \newtheorem{exercise}{Exercise}[section]
\fi
\ifx\conjecture\undefined
  \newtheorem{conjecture}{Conjecture}[section]
\fi

\theoremstyle{remark}
\ifx\remark\undefined
  \newtheorem{remark}{Remark}
\fi

\usepackage[capitalise,nameinlink,noabbrev]{cleveref}
\crefformat{equation}{Eq.\ (#2#1#3)}
\crefrangeformat{equation}{Eqs.\ (#3#1#4) to~(#5#2#6)}
\crefmultiformat{equation}{Eq.\ (#2#1#3)}%
{ and~(#2#1#3)}{, (#2#1#3)}{ and~(#2#1#3)}
\crefname{assumption}{Assumption}{Assumptions}
\crefname{theorem}{Theorem}{Theorems}
\crefname{proposition}{Proposition}{Propositions}
\crefname{corollary}{Corollary}{Corollaries}
\crefname{lemma}{Lemma}{Lemmas}
\crefname{remark}{Remark}{Remarks}
\crefname{definition}{Definition}{Definitions}
\crefname{example}{Example}{Examples}

\newcommand{\reels}{\mathbb{R}}
\newcommand{\naturels}{\mathbb{N}}

\newcommand{\proba}{\mathbb{P}}
\newcommand{\Tau}{\mathrm{T}}

\endlocaldefs

\begin{document}

\begin{frontmatter}
\title{Approximation convergence in the inverse first-passage time problem}
\runtitle{Approximation in the inverse first-passage time problem}

\begin{aug}

\author[A]{\fnms{Yoann}~\snm{Potiron}\ead[label=e1]{potiron@fbc.keio.ac.jp}},
\address[A]{Faculty of Business and Commerce, Keio University\printead[presep={,\ }]{e1}}
\end{aug}

\begin{abstract}
The inverse first-passage time problem determines a boundary such that the first-passage time of a Wiener process to this boundary has a given distribution. An approximation which is based on the starting value of the boundary to a smooth boundary by a piecewise linear boundary is given by equating the probability of the first-passage time to a linear boundary and the increment of the distribution on each interval. We propose a modification of that approximation which also approximates the starting value of the boundary. First, we show that the approximation is well-defined when assuming that the boundary is absolutely continuous. Second, we show that a subsequence of this new approximation uniformly converges to the boundary when the length of each interval of linear approximation goes to 0 asymptotically. The results are obtained using Arzel\`a-Ascoli theorem on any compact space on which we further assume that the boundary admits uniformly dominated derivative. As the starting value of the boundary is unknown, this makes the new approximation more suitable for applications. The results are also proved in the first-passage time problem of a reflected Wiener process. 

\end{abstract}



\end{frontmatter}

\section{Introduction}
The first-passage time (FPT) problem in statistics can at least be traced back to the one-sample Kolmogorov-Smirnov statistic for which the stochastic process is equal to the difference between the true and empirical cumulative distribution function (cdf). Explicit solutions of the distribution are only found for the linear boundary (see \cite{doob1949heuristic} (Equation (4.2), p. 397) or \cite{malmquist1954certain} (p. 526)), the upper and lower linear boundary (see \cite{anderson1960modification} (Theorem 5.1, p. 191)), the Daniels boundary (see \cite{daniels1969minimum}), a general boundary setup but which depends on asymptotic conditional expectations (see \cite{durbin1985first}), the one-sided square root boundary (see \cite{novikov1971stopping}), the quadratic boundary (see \cite{salminen1988first}), the piecewise-linear boundary (see \cite{wang1997boundary}) and the piecewise-specific boundary (\cite{novikov1999approximations}) where "specific" can refer to any of the aforementioned cases. \cite{mehr1965certain} consider a Gauss-Markov stochastic process. \cite{lai1974control}, \cite{lai1977nonlinear}, \cite{lai1979nonlinear}, \cite{gut1974moments}, \cite{woodroofe1976renewal} and \cite{woodroofe1977second} consider the case when the stochastic process is a discrete-time sum of i.i.d. variables. \cite{siegmund1986boundary} develops tools to calculate FPT cdf when the stochastic process is a Wiener process and discusses applications in sequential analysis. \cite{matthews1985asymptotic} show that tests based on maximal score statistics involve the solution to a FPT problem of an Ornstein-Uhlenbeck process. In survival analysis, \cite{butler1997stochastic} present a Bayesian approach when the stochastic process is a semi-Markov process. In econometrics, \cite{abbring2012mixed} and \cite{renault2014dynamic} study mixed FPT where the stochastic process is respectively a spectrally negative Levy process and the sum of a Wiener process and a positive linear trend.

The inverse first-passage problem determines the boundary function such that the first-passage time of a standard Wiener process to this boundary has a given distribution. Almost 50 years ago, A. Shiryaev, during a Banach center meeting in 1976, asked if one can determine a continuous boundary with exponential distribution, which is commonly referred as the inverse Shiryaev problem. \cite{dudley1977stopping} show the existence of a stopping time with respect to a general stochastic
process, but this stopping time is not a FPT. The existence of lower semi-continuous solutions was established in \cite{anulova1981markov} for the FPT of a reflected Wiener process by compacity arguments in a
discrete approximation of the boundary and the distribution. When the distribution is non-atomic, \cite{cheng2006analysis} and \cite{chen2011existence} show the existence and uniqueness of the inverse first-passage problem for diffusions by a transfer into a free boundary problem. For a general distribution, \cite{ekstrom2016inverse} show the existence and uniqueness for Wiener processes by discretizing an optimal stopping problem. \cite{beiglbock2018geometry} consider a more general optimal stopping problem which yields existence and uniqueness as a by-product. \cite{chen2022higher} study higher-order regularity properties
of the solution of the inverse first-passage problem. The uniqueness for reflected Wiener processes is shown by a discrete approximation argument along with stochastic ordering in \cite{klump2023uniqueness}. The existence and the uniqueness for Levy processes and diffusions are studied in \cite{klump2023conditions}. There are other related papers. The problem was reformulated as a nonlinear Volterra integral equation in \cite{peskir2002integral}. \cite{peskir2002limit} study the behavior in the neighborhood of $0$. \cite{abundo2006limit} consider extensions to the general diffusion process case.

As explicit solutions are only found in a few cases, the literature related to the inverse first-passage problem relies heavily on approximation as in \cite{zucca2009inverse} and \cite{song2011approximation}. Based on \cite{wang1997boundary}, an approximation to a continuous boundary by a piecewise linear boundary is given by equating the probability of the FPT to a linear boundary and the increment of the cdf on each interval in \cite{zucca2009inverse}. That approximation is based on the starting value of the boundary, which has to be guessed in practice since it is unknown. Assuming that the FPT cdf is absolutely continuous and the boundary is concave or convex, they prove that the error due to the approximation does not explode but that it is dominated by the initial error. We propose a modification of that approximation which also approximates the starting value of the boundary, which makes it more suitable for applications. First, we show that the new approximation is well-defined when the boundary is absolutely continuous. This is a rigorous result of Remark 3.2 in \cite{zucca2009inverse} which also shows that the starting value of the new approximation is well-defined. Also, we have more explicit assumptions than the cited paper who assumes that all regularity assumptions ensuring the existence of the objects introduced and properties imposed are fulfilled.  Second, we show that a subsequence of this new approximation uniformly converges to the boundary when the length of each interval of linear approximation goes to 0 asymptotically. The results are obtained using Arzel\`a-Ascoli theorem on any compact space on which we assume that the boundary is absolutely continuous with uniformly dominated derivative. The use of an asymptotic when the length of each interval of linear approximation goes to 0 and the convergence result are new to the literature on inverse FPT problem and important as they indicate that this new approximation can be relatively safely used to approximate the unknown boundary when choosing a small length in practice. The results are also proved in the FPT problem of a reflected Wiener process, which are also new.



\section{Setting and main results}
We consider the complete stochastic basis $\mathcal{B} = (\Omega, \proba, \mathcal{F}, \mathbf{F})$, where $\mathcal{F}$ is a $\sigma$-field and $\mathbf{F} = (\mathcal{F}_t)_{t \in \reels^+}$ is a filtration. For $t_0 \geq 0$, $A \subset \reels^+$ and $B \subset \reels$ such that $t_0 \in A$, we define the set of continuous functions with non-negative starting values at time $t_0$ as $\mathcal{C}_{t_0}^+(A,B)  =  \{h: A \rightarrow B  \text{ s.t. } h \text{ is continuous and } h(t_0) \geq 0\}$. We first give the definition of the set of boundary functions. Since the approximation by a piecewise linear boundary given in \cite{wang1997boundary} requires continuity of the boundary, we restrict ourselves to the continuous boundary case.

\begin{definition}\label{defboundaryset}
We define the set of boundary functions started at time $t_0 \geq 0$ as $\mathcal{G} = \mathcal{C}_{t_0}^+(\reels^+,\reels).$
\end{definition}
\begin{definition}\label{defFPT}
We define the FPT of an $\mathbf{F}$-adapted continuous process $(Z_t)_{t \in \reels^+}$ started at time $t_0 \geq 0$ to a boundary
$g \in \mathcal{G}$ as
\begin{eqnarray}
\label{TgZdef}
\Tau_g^Z = \inf \{t \in \reels^+ \text{ s.t. } Z_t - Z_{t_0} \geq g(t_0+t)\}.
\end{eqnarray}
\end{definition}
\noindent We define an $\mathbf{F}$-standard Wiener process as $(W_t)_{t \in \reels^+}$. We will consider the two cases in the following of this paper: 
\begin{enumerate}
    \item (Wiener process) $Z_t=W_t$
    \item (reflected Wiener process) $Z_t=|W_t|$
\end{enumerate}
Since the stochastic process $Z$ is continuous and $\mathbf{F}$-adapted and the set generated by the points above the boundary is optional, the first-hitting time $\Tau_{g}^{Z}$ is a $\mathbf{F}$-stopping time by Theorem I.1.27 (p. 7) in \cite{JacodLimit2003}. We define the cdf of $Z$ as
\begin{eqnarray}
\label{PgZdef}
P_g^Z(t)= \proba (\Tau^Z_g \leq t) \text{ for any } t \geq 0.
\end{eqnarray}
The basic assumption for the approximation by a piecewise linear boundary given in \cite{zucca2009inverse} is that $P_g^Z$ is absolutely continuous. Accordingly, the authors assume that all regularity assumptions ensuring the existence of the objects introduced and properties imposed are fulfilled. When $g$ is continuous, we know by Theorem 8.1 in \cite{ekstrom2016inverse} that $P_g^Z$ is continuous. When $g$ is continuous with continuous derivative, we know by Lemma 3.3 in \cite{strassen1967almost} that $P_g^Z$ is continuous with continuous derivative. In fact, we can show that when $g$ is absolutely continuous, then $P_g^Z$ is absolutely continuous. 

\begin{proof}[\textbf{Assumption A}]
We assume that $g$ is absolutely continuous on $[t_0,\infty)$ and that $\underset{t \rightarrow t_0, t>t_0}{\liminf} g'(t) > 0$ if $g(t_0)=0$.
\phantom\qedhere
\end{proof}

\begin{proposition}
\label{propstrassen}
We assume that \textbf{Assumption A} holds. Then, $P_g^Z$ is absolutely continuous on $\reels^+$.
\end{proposition}
\noindent Since $P_g^Z(t)$ is absolutely continuous, there exists a pdf $f_g^Z:\reels^+ \rightarrow \reels^+,$ defined as
\begin{eqnarray}
\label{fZgt}
f_g^Z(t)= \frac{dP_g^Z(t)}{dt} \text{ for any } t \geq 0.
\end{eqnarray}
By Proposition \ref{propstrassen} along with \textbf{Assumption A}, there is no loss of generality in restricting ourselves to the absolute continuous cdf case.
\begin{definition} \label{Fdef}
A function $F:\reels^+ \rightarrow  [0,1]$ is a survival cdf if $F$ is nondecreasing, absolutely continuous, i.e., with pdf $f:\reels^+ \rightarrow \reels^+,$ defined as
\begin{eqnarray}
\label{f}
f(t)= \frac{dF(t)}{dt} \text{ for any } t \geq 0,
\end{eqnarray} 
and satisfies $F(0) = 0$ and $\underset{t \rightarrow \infty}{\lim} F(t) =1$. 
\end{definition}
\noindent The inverse first-passage time problem determines a boundary $g \in \mathcal{G}$ such that
\begin{eqnarray}
\label{ffg}
f_g^Z(t)= f(t) \text{ for any } t \geq 0.
\end{eqnarray}
As explicit solutions are only found in a few cases, the literature related to the inverse first-passage problem relies heavily on approximations. Based on \cite{wang1997boundary} idea, an approximation to a continuous boundary by a piecewise linear boundary is given in \cite{zucca2009inverse}. More specifically, they consider a time discretization $t_m = t_0 + m \Delta$, where $t_0 \geq 0$ is the starting time, $\Delta > 0$ is the length of each interval of linear approximation and $m \in \naturels$. That approximation is based on the starting value of the boundary, i.e., $g(0)$. Since $g(0)$ is unknown, it has to be guessed in practice. The driving idea of their algorithmic approximation is to determine recursively the slope $\alpha_m$ of the linear increment approximation on $[t_{m}, t_{m+1}]$ by equating the probability of the
 FPT of $W$ to the approximation and the increment of the survival cdf, i.e., $\int_{t_m}^{t_{m+1}}f(s)ds$. Although their approximation works for any fixed $\Delta > 0$, they do not propose any asymptotic result when $\Delta \rightarrow 0$. First, we propose a modification of that approximation which also approximates the starting value of the boundary, which makes it more suitable for applications. Second, we consider an asymptotics where the length of each interval of linear approximation goes to 0, i.e., $\Delta \rightarrow 0$. We define $t_f \in \overline{\reels}_*^+$ as the final time where $\overline{\reels}_*^+ = \reels_*^+ \cup \infty$ and $t_0 < t_f$. For any $n \in \naturels_*$, we consider a time discretization $t_m^n = t_0^n + m \Delta_n$ for any $m \in \naturels$ such that $t_m^n\leq t_f$ where $t_0^n=t_0$ and $\Delta_n = \frac{(t_f-t_0)\mathbf{1}_{\{t_f<\infty \}}+\mathbf{1}_{\{t_f=\infty \}}}{2^n}$. We consider a discretization length in the order $\frac{1}{2^n}$ so that we obtain that the time discretization is nested, i.e., for any $t_m^n$ and any $l \geq m$ there exists a time $t_k^l$ such that $t_m^n=t_k^l$. This is required in our proofs to prove that the limit of a subsequence obtained by Arzel\`a-Ascoli theorem satisfies Equation (\ref{ffg}).
 \begin{definition}\label{defpiecewiselinearboundaryset}
For any $n \in \naturels_*$, we define the subset of piecewise linear boundary functions as 
$$\mathcal{G}_{PL}^{n} = \Big\{g \in \mathcal{G} \text{ s.t. } g \text{ is linear on each interval } \big[t_m^n, t_{m+1}^n\big] \Big\}.$$
\end{definition}
\noindent For any $n \in \naturels_*$, we define the sequence of piecewise linear approximation of the boundary $g^{n} \in \mathcal{G}_{PL}^{n}$ recursively on $m$ as
\begin{eqnarray}
\label{seq1}
g^{n}(0) & = & \alpha_{0}^n,\\
\label{seq2}
g^{n}(u) & = & g^{n}(t_m^n) + \alpha_{m+1}^n (u - t_m^n) \text{ for any } u \in (t_m^n,t_{m+1}^n], m \in \naturels \text{ s.t. } t_{m+1}^n\leq t_f,
\end{eqnarray}
with $\alpha_{0}^n \in \reels_*^+$ and $\alpha_{m}^n \in \reels$ for $m \neq 0$ satisfying 
\begin{eqnarray}
\label{alpha0n} \proba \big(T_{\alpha_{0}^n}^Z \in [0,\Delta_n]\big) &= &\int_{0}^{\Delta_n}f(s)ds,\\
\label{pintf}\proba \big(T_{g^{n}}^Z \in [m\Delta_n,(m+1)\Delta_n]\big) & = & \int_{m\Delta_n}^{(m+1)\Delta_n} f(s)ds \text{ for any } m \in \naturels \text{ s.t. } t_{m+1}^n\leq t_f.
\end{eqnarray}
Equations (\ref{seq1})-(\ref{seq2}) and Equation (\ref{pintf}) correspond exactly to Equations (3.1)-(3.2) in \cite{zucca2009inverse}. The novelty in this paper is Equation (\ref{alpha0n}) in which we determine the approximation of the starting value of the boundary, i.e., $\alpha_{0}^n$, by equating the probability of the
 FPT of $Z$ to a constant boundary equal to $\alpha_{0}^n$ and the increment of the survival cdf on the first interval. The next lemma gives a more explicit form to $\alpha_{0}^n$ since the only unknown value in Equation (\ref{alpha0neq2sim}) is $\alpha_{0}^n$. We define the probability of the FPT started at time $t_0^n$ to a constant boundary equal to $\alpha \in \reels^+_*$ on $[t_0^n,t_1^n]$ as $G_{0}^n : \reels^+_* \rightarrow \reels$ such that
\begin{eqnarray}
\label{defF0n}
G_{0}^n(\alpha) = 1 - \int_{-\infty}^{\alpha} \Big(1 - \exp\Big(\frac{-2(\alpha-x_{1})\alpha}{\Delta_n}\Big)\Big) \frac{\exp (-\frac{x_{1}^2}{\Delta_n})}{\sqrt{\pi \Delta_n}} dx_1.
\end{eqnarray}
\begin{lemma}
\label{lemmaexpform0}
For any $n \in \naturels_*$, Equation (\ref{alpha0n}) can be reexpressed as
\begin{eqnarray}
\label{alpha0neq2sim}
G_{0}^n(\alpha_{0}^n)  - \int_{0}^{\Delta_n}f(s)ds = 0.
\end{eqnarray}
\end{lemma}
\noindent In the next lemma, we give a more explicit form to $\alpha_1^n$ from Equation (\ref{pintf}) based on the known value $\alpha_0^n$. We define the probability of the FPT started at time $t_0^n$ from $\alpha_{0}^n$ to a linear boundary with trend $\alpha \in \reels$ on $[t_0^n,t_1^n]$ as $G_{1}^n : \reels^+ \rightarrow \reels$ such that
\begin{eqnarray}
\label{defF1n}
G_{1}^n(\alpha) = 1 - \int_{-\infty}^{\alpha_0^n + \alpha \Delta_n} \Big(1 - \exp\Big(\frac{-2(\alpha_0^n + \alpha \Delta_n-x_{1})\alpha_0^n}{\Delta_n}\Big)\Big) \frac{\exp (-\frac{x_{1}^2}{\Delta_n})}{\sqrt{\pi \Delta_n}} dx_1.
\end{eqnarray}
\begin{lemma}
\label{lemmaexpform1}
For any $n \in \naturels_*$ and $m=1$, Equation (\ref{pintf}) can be reexpressed as
\begin{eqnarray}
\label{alpha01neq2sim}
G_{1}^n(\alpha_{1}^n)  - \int_{0}^{\Delta_n}f(s)ds = 0.
\end{eqnarray}
\end{lemma}
\noindent In the next lemma, we give a more explicit form to $\alpha_2^n$ from Equation (\ref{pintf}) based on the known values $\alpha_{1}^n$ and $\alpha_{0}^n$. We define the probability of the FPT started at time $t_1^n$ from $g^n(t_1^n)$ to a linear boundary with trend $\alpha \in \reels$ on $[t_1^n,t_2^n]$ as $G_{2}^n : \reels^+ \rightarrow \reels$ such that
\begin{eqnarray}
\label{defF2n}
G_{2}^n(\alpha) & = & 1 - \int_{-\infty}^{\alpha_0^n+\alpha_1^n \Delta_n} \int_{-\infty}^{\alpha_0^n+(\alpha_1^n + \alpha) \Delta_n} \Big(1 - \exp\Big(\\ \nonumber & & \frac{-2(\alpha_0^n+(\alpha_1^n + \alpha) \Delta_n-x_{2})(\alpha_0^n+\alpha_1^n \Delta_n-x_{1})}{\Delta_n}\Big)\Big) \frac{\exp (-\frac{(x_{2}-x_{1})^2}{\Delta_n})}{\sqrt{\pi\Delta_n}}\\ \nonumber & &   
\Big(1 - \exp\Big(\frac{-2(\alpha_0^n+\alpha_1^n \Delta_n-x_{1})\alpha_0^n}{\Delta_n}\Big)\Big) \frac{\exp (-\frac{x_{1}^2}{\Delta_n})}{\sqrt{\pi\Delta_n}} dx_1 dx_2.
\end{eqnarray}
\begin{lemma}
\label{lemmaa1n}
For any $n \in \naturels_*$ and $m=2$, Equation (\ref{pintf}) can be reexpressed as
\begin{eqnarray} \nonumber
\label{alpha1neq2sim} G_{2}^n(\alpha_{2}^n)  - \int_{0}^{2 \Delta_n}f(s)ds = 0.
\end{eqnarray}
\end{lemma}
\noindent In the lemma that follows, we give a more explicit form to $\alpha_m^n$ from Equation (\ref{pintf}) in the case $m>2$ based on known values $(\alpha_{k}^n)_{k=0, \ldots, m-1}$. For any $m \in \naturels$ such that $m>2$ and $t_{m+1}^n\leq t_f$, we define $x_0=0$ and the probability of the FPT started at time $t_{m-1}^n$ from $g^n(t_{m-1}^n)$ to a linear boundary with trend $\alpha \in \reels$ on $[t_{m-1}^n,t_m^n]$ as $G_{m}^n : \reels^+ \rightarrow \reels$ such that
\begin{eqnarray}
\nonumber
G_{m}^n(\alpha) & = & 1 - \int_{-\infty}^{\alpha_0^n+\alpha_1^n \Delta_n} \int_{-\infty}^{\alpha_0^n+(\alpha_1^n + \alpha_2^n) \Delta_n} \ldots \int_{-\infty}^{\alpha_0^n+(\sum_{k=1}^{m}\alpha_k^n + \alpha) \Delta_n} (1 - \\ \nonumber & & \exp(\frac{-2(\alpha_0^n+(\sum_{i=1}^{m}\alpha_i^n + \alpha) \Delta_n-x_{m+1})(\alpha_0^n+\sum_{i=1}^{m}\alpha_i^n \Delta_n-x_{m})}{\Delta_n})) \\ \nonumber  & &\frac{\exp (-\frac{(x_{m+1}-x_{m})^2}{\Delta_n})}{\sqrt{\pi\Delta_n}} \prod_{k=0}^{m-1} (1 - \\ \nonumber & & \exp(\frac{-2(\alpha_0^n+\sum_{i=1}^{k+1}\alpha_i^n \Delta_n-x_{k+1})(\alpha_0^n+\sum_{i=1}^{k}\alpha_i^n \Delta_n-x_{k})}{\Delta_n})) \\ \label{defFmn}  & &\frac{\exp (-\frac{(x_{k+1}-x_{k})^2}{\Delta_n})}{\sqrt{\pi\Delta_n}} dx_1 dx_2\ldots dx_{m+1}.
\end{eqnarray}
\begin{lemma}
\label{lemmaexpformm}
For any $n \in \naturels_*$ and any $m \in \naturels$ such that $m>2$ and $t_{m+1}^n\leq t_f$, Equation (\ref{pintf}) can be reexpressed as
\begin{eqnarray} G_{m}^n(\alpha_{m}^n)  - \int_{0}^{m \Delta_n}f(s)ds = 0. \label{alphamneq2sim}
\end{eqnarray}
\end{lemma}
\noindent With the same arguments as in the proofs of Lemmas \ref{lemmaexpform0}-\ref{lemmaexpform1}-\ref{lemmaa1n}-\ref{lemmaexpformm}, for any $n \in \naturels_*$ we can define $H_{0}^n$ as $H_{0}^n : \reels^+_* \rightarrow \reels$
 and  $H_{m}^n$ as $H_{m}^n : \reels^+ \rightarrow \reels$ for any $m \in \naturels_*$ such that $t_{m}^n\leq t_f$ which are defined in the same way as $G_{m}^n$ for the reflected Wiener process. As the obtained equations are much longer than in the Wiener process case, we do not report them. Our next result establishes that the sequence $\alpha_{m}^n$ is well-defined. This is a rigorous result of Remark 3.2 in \cite{zucca2009inverse}. This also includes that $\alpha_{0}^n$ is well-defined, which is not treated in the cited paper. Moreover, we have more explicit assumptions than the cited paper who assumes that all regularity assumptions ensuring the existence of the objects introduced and properties imposed are fulfilled. The main idea of the proof is based on Equations (3.4)-(3.5) from the cited paper. 
\begin{proposition}
\label{propexistence}
We assume that \textbf{Assumption [A]} holds. For any $n \in \naturels^*$, Equation (\ref{alpha0n}) defines a unique $\alpha_{0}^{n} \in \reels_*^+$ and Equation (\ref{pintf}) defines a unique $\alpha_{m}^{n} \in \reels$ for any $m \in \naturels$ such that $t_m^n\leq t_f$. 
\end{proposition}
\noindent We give our main result in the next theorem. This shows that a subsequence of the approximation uniformly converges to the boundary when the length of each interval of linear approximation goes to 0 asymptotically. The results are obtained using Arzel\`a-Ascoli theorem on any compact space $[t_0,t_f]$, so we make the following assumption. 
\begin{proof}[\textbf{Assumption B}]
We assume that $t_f < \infty$.  \phantom\qedhere
\end{proof}
\noindent Finally, we assume that the boundary has a uniformly dominated derivative on $[t_0,t_f]$. This allows us to show that the $\alpha_{m}^{n}$ are uniformly bounded, which in turn implies that the approximated boundary is uniformly bounded and uniformly equicontinuous.
\begin{proof}[\textbf{Assumption C}]
We assume that $g'(t)$ exists and is uniformly dominated on $[t_0,t_f]$, i.e. that $\underset{t \in [t_0,t_f]}{\sup} \mid g'(t) \mid < \infty$.
\phantom\qedhere
\end{proof}
\begin{theorem}
\label{main}
We assume that \textbf{Assumption [A]}, \textbf{Assumption [B]} and \textbf{Assumption [C]} hold. Then, there exists a subsequence $g^{n_k}$ of $g^n$ which converges uniformly to $g$ on $[t_0,t_f]$, i.e., $\underset{t \in [t_0,t_f]}{\sup} \mid g^{n_k}(t) - g(t) \mid \rightarrow 0$ as $n\rightarrow \infty$.
\end{theorem}

\section{Proof that the approximation is well-defined}
In this section, we prove that $\alpha_{m}^n$ is well-defined. We start with the proof of Proposition \ref{propstrassen}, which extends the arguments from the proof of Lemma 3.3 in \cite{strassen1967almost}.
\begin{proof}[Proof of Proposition \ref{propstrassen}] By \textbf{Assumption [A]}, we have that $g$ is absolutely continuous on $[t_0,\infty)$. Thus $g$ admits a derivative almost everywhere on $[t_0,\infty)$, i.e., there exists a Lebesgue-negligible set $\mathcal{N} \subset [t_0,\infty)$ such that $g$ admits a derivative for any $t \in [t_0,\infty) - \mathcal{N}$. we define the set of linear transformation from $\mathcal{N}$ as $\widetilde{\mathcal{N}} \subset \reels^+$ such that $t \in \widetilde{\mathcal{N}}$ if $t+t_0 \in \mathcal{N}$. Since $\mathcal{N}$ is a Lebesgue-negligible set, we have by construction that $\widetilde{\mathcal{N}}$ is a Lebesgue-negligible set. We show now that $P_g^W$ is absolutely continuous on $\reels^+$. Since $\widetilde{\mathcal{N}}$ is a Lebesgue-negligible set, it is sufficient to show that $P_g^W$ admits a derivative for any $t \in \reels^+ - \widetilde{\mathcal{N}}$. It is then sufficient to show that $P_g^W$ admits a derivative on any open interval $(u,v)$ where $u \in \reels^+$ and $v \in \reels^+$ satisfy $u < v$ and $(u,v) \cap \widetilde{\mathcal{N}} = \emptyset$. We can show this statement by extending the arguments from the proof of Lemma 3.3 in \cite{strassen1967almost} along with the assumption that $\underset{t \rightarrow t_0, t>t_0}{\liminf} g'(t) > 0$ if $g(t_0)=0$ by \textbf{Assumption [A]}. The reflected Wiener process case follows since the FPT of a reflected Wiener process to a linear boundary is equal to the FPT of a Wiener process to a symmetric upper linear boundary and lower linear boundary when the boundary from the reflected Wiener process and the upper boundary are equal.

\end{proof}
\begin{definition}
We define the transition pdf of the stochastic process $Z$ at
time $t$ constrained by the absorbing boundary $g$ over $[s, t]$ given that $Z_s = y$ as $p_g^Z(t, x \mid s, y)$ such that
\begin{eqnarray}
\label{constrainedtransitionproba}
p_g^Z(t, x \mid s, y) = \frac{\partial}{\partial x}\proba(Z_t < x, \Tau_g^Z > t \mid Z_s = y)
\end{eqnarray}
with $x < g(t)$, $t > s \geq 0$ and $y < g(s)$ given and fixed. 
\end{definition}
\noindent In the following lemma, we give the pdf and the transition pdf for the FPT of a Wiener process to a linear boundary. This is a consequence to \cite{doob1949heuristic} (Equation (4.2), p. 397), \cite{malmquist1954certain} (p. 526) and \cite{durbin1971boundary} (Lemma 1).
\begin{lemma} \label{exlinear}
We assume that 
$$g(t)= \alpha_1 t + \alpha_0 \text{ for any } t \geq 0,$$
where $t_0 \geq 0$, $\alpha_0 \in \reels_*^+$, $\alpha_1 \in \reels$, $x_0 \in \reels$ such that $g(t_0) > x_0$. We have that the pdf is equal to
\begin{eqnarray}
\label{densitylinearupperboundary}
f_g^W(t \mid t_0, x_0) & = &  \frac{\alpha_0-x_0}{\sqrt{2\pi (t-t_0)^3}}\exp \Big(-\frac{(\alpha_0+\alpha_1(t-t_0)-x_0)^2}{2(t-t_0)} \Big).
\end{eqnarray}
The transition pdf is equal to
\begin{eqnarray}
\label{exlineareq} p_g^W(t_1, x_1 \mid t_0, x_0) & = & \Big(1 - \exp \Big(\frac{-2(g(t_1)-x_1)(g(t_0)-x_{0})}{t_1-t_{0}}\Big)\Big) \frac{\exp (-\frac{(x_1-x_{0})^2}{2(t_1-t_{0})})}{\sqrt{2\pi(t_1-t_{0})}}.
\end{eqnarray}
\end{lemma}
\begin{proof}[Proof of Lemma \ref{exlinear}]
Equation (\ref{densitylinearupperboundary}) is obtained in \cite{doob1949heuristic} (Equation (4.2), p. 397) or \cite{malmquist1954certain} (p. 526). Equation (\ref{exlineareq}) follows from \cite{durbin1971boundary} (Lemma 1).
\end{proof}
\noindent In the following lemma, we give the pdf and the transition pdf for the FPT of a Wiener process to a continuous piecewise linear boundary. Equation (\ref{exeqpiecewiselinear}) is already available in \cite{wang1997boundary} and \cite{zucca2009inverse} (Section 2.1.3, pp. 1323-1324)
\begin{lemma} \label{expiecewiselinear} We assume that 
$$g(t)= \alpha_i t + \beta_i, \text{ for any } t \in [t_{i-1}, t_i] $$
with $t_i =  i\Delta + t_0$, where $t_0 \geq 0$, $\Delta > 0$ and $\alpha_i, \beta_i \in \reels$ satisfying $\alpha_{i+1} + \beta_{i+1} t_i = \alpha_i + \beta_i t_i$ so that the boundary is continuous. We can reexpress the transition pdf as
\begin{eqnarray}
\label{pgWfact} p_g^W(t_1, x_1, \ldots, t_n, x_n \mid t_0, x_0) & = & \prod_{i=1}^{n} p_g^W(t_i, x_i \mid t_{i-1}, x_{i-1}),
\end{eqnarray}
with $(x_1, x_2, \ldots, x_n) \in \reels^n$ and $x_i \leq g(t_i)$ for $i=1,\ldots,n$ and $x_0 < g(t_0)$ where $t_0 < t_1 < t_2 < \ldots < t_n$ are given and fixed. We can reexpress the transition pdf with an explicit expression as
\begin{eqnarray}
\label{exeqpiecewiselinear}
p_g^W(t_1, x_1, \ldots, t_n, x_n \mid t_0, x_0) & = & \prod_{i=1}^{n} (1 - \exp((\frac{-2(g(t_i)-x_i)(g(t_{i-1})-x_{i-1})}{t_i-t_{i-1}}))) \\ \nonumber & & \frac{\exp (-\frac{(x_i-x_{i-1})^2}{2(t_i-t_{i-1})})}{\sqrt{2\pi(t_i-t_{i-1})}}.
\end{eqnarray}
We can deduce that
\begin{eqnarray}
\label{exeqpiecewiselinear2}
\proba(W_{t_1} \in C_1, \ldots, W_{t_n} \in C_n, \Tau_g^W > t_n \mid W_{t_0} = x_0)\\
\nonumber  =  \int_{C_1} \ldots \int_{C_n} p_g^W(t_1, x_1, \ldots, t_n, x_n \mid t_0, x_0) dx_1 \ldots dx_n
\end{eqnarray}
for any Borel set $C_i \subset (-\infty, g(t_i)]$ with $i=1,\ldots,n$.
\end{lemma}
\begin{proof}[Proof of Lemma \ref{expiecewiselinear}] Equation (\ref{pgWfact}) is obtained by Definition (\ref{constrainedtransitionproba}) and follows by induction with conditional probability formula. Then, Equation (\ref{exeqpiecewiselinear}) can be deduced by plugging Equation (\ref{exlineareq}) into Equation (\ref{pgWfact}). Finally, Equation (\ref{exeqpiecewiselinear2}) is a direct consequence of Equation (\ref{exeqpiecewiselinear}).
\end{proof}
\noindent In the following lemma, we give the pdf for the FPT of a reflected Wiener process to a linear boundary. This is based on the explicit solution from \cite{anderson1960modification} (Theorem 5.1, p. 191) for the FPT to an upper linear boundary and a lower linear boundary. This is due to the fact that the FPT of a reflected Wiener process to a linear boundary is equal to the FPT of a Wiener process to a symmetric upper linear boundary and lower linear boundary when the boundary from the reflected Wiener process and the upper boundary are equal. Note that we could deduce the transition pdf and transition pdf for the piecewise linear boundary with the same arguments as for the proofs of Lemma \ref{exlinear} and Lemma \ref{expiecewiselinear}.
\begin{lemma} \label{exlinearsym}
We assume that 
$$g(t)= \alpha_1 t + \alpha_0 \text{ for any } t \geq 0,$$
where $t_0 \geq 0$, $\alpha_0 \in \reels_*^+$, $\alpha_1 \in \reels$, $x_0 \in \reels$ such that $g(t_0) > x_0$. We have that the pdf is equal to
\begin{eqnarray}
\label{fgWex20symcond} f_g^{ \mid W \mid} (t_0 \mid t_0, x_0) & = & 0,\\
\label{fgWex21symcond} f_g^{ \mid W \mid} (t \mid t_0, x_0) & = & \frac{2}{(t-t_0)^{3/2}} \phi (\frac{\alpha_1 t + \alpha_0-x_0}{\sqrt{t-t_0}})\sum_{r=0}^{\infty} \Big\{ (4 r +1) (\alpha_0-x_0)  \\ \nonumber & & e^{-(8r(r+1) (\alpha_0-x_0)) (\alpha_1 t+\alpha_0-x_0) /t)} - (4r +2) (\alpha_0-x_0)\\ \nonumber & &  e^{-(4(r+1)(2r+1) (\alpha_0-x_0)(\alpha_1(t-t_0) + \alpha_0-x_0) /(t-t_0))}
\Big\} \text{ for any } t>t_0,  
\end{eqnarray}
where $\phi$ is defined as the standard Gaussian density function.
\end{lemma}

\begin{proof}[Proof of Lemma \ref{exlinearsym}] We first consider the FPT of a Wiener process to an upper linear boundary and a lower linear boundary. We first assume that the boundary is upper linear and lower linear, i.e., that
\begin{eqnarray*}
g(t) & = & (\gamma_2 + \delta_2 t, \gamma_1 + \delta_1 t),
\end{eqnarray*}
where $\gamma_1 > 0$, $\gamma_2 < 0$, $\delta_1 \geq \delta_2$ and not $\delta_1 = \delta_2 = 0$. By \cite{anderson1960modification} (Theorem 5.1, p. 191), we have that the FPT pdf is equal to 
\begin{eqnarray}
\label{fgWex20} f_g^W (0) & = & 0,\\
\label{fgWex21} f_g^W (t) & = & \frac{1}{t^{3/2}}\Bigg[ \phi (\frac{\delta_1 t + \gamma_1}{\sqrt{t}})\sum_{r=0}^{\infty} \Big\{ ((2 r +1) \gamma_1 - 2r \gamma_2) \\ \nonumber & & e^{-(2r /t)(r \gamma_1 - (r+1) \gamma_2)(\delta_1 t + \gamma_1-(\delta_2 t + \gamma_2))}\\ \nonumber & & - (2 (r +1) \gamma_1 - 2r \gamma_2) e^{-(2(r+1) /t)((r+1) \gamma_1 - r \gamma_2)(\delta_1 t + \gamma_1-(\delta_2 t + \gamma_2))}
\Big\} \\ \nonumber
& & + \phi (\frac{\delta_2 t + \gamma_2}{\sqrt{t}})\sum_{r=0}^{\infty} \Big\{ (2 r \gamma_1 - (2r + 1) \gamma_2) \\ \nonumber & & e^{-(2(r+1) /t)((r+1) \gamma_1 - r \gamma_2)(\delta_1 t + \gamma_1-(\delta_2 t + \gamma_2))} \\ \nonumber & & - (2 (r +1) \gamma_1 - 2r \gamma_2) e^{-(2r /t)(r \gamma_1 - (r+1) \gamma_2)(\delta_1 t + \gamma_1-(\delta_2 t + \gamma_2))}
\Big\} \Bigg]  
\end{eqnarray}
for any $t>0$. Now we assume that the boundaries are symmetric, i.e., that $g= (-\alpha_1 t - \alpha_0, \alpha_1 t + \alpha_0)$ where $\alpha_1 \in \reels$ and $D \in \reels^+_*$. From Equations (\ref{fgWex20})-(\ref{fgWex21}), we can deduce that 
\begin{eqnarray}
\label{fgWex20sym} f_g^W (0) & = & 0,\\
\label{fgWex21sym} f_g^W (t) & = & \frac{2}{t^{3/2}} \phi (\frac{\alpha_1 t + \alpha_0}{\sqrt{t}})\sum_{r=0}^{\infty} \Big\{ (4 r +1) \alpha_0  e^{-(8r(r+1) \alpha_0) (\alpha_1 t+\alpha_0) /t)}\\ \nonumber & & - (4r +2) \alpha_0 e^{-(4(r+1)(2r+1) \alpha_0(\alpha_1 t + \alpha_0) /t)}
\Big\} \text{ for any } t>0.  
\end{eqnarray}
From Equations (\ref{fgWex20sym})-(\ref{fgWex21sym}) and since the FPT of a reflected Wiener to a linear boundary is equal to the FPT of a Wiener process to a symmetric upper linear boundary and lower linear boundary when the boundary from the reflected Wiener process and the upper boundary are equal, we can deduce Equations (\ref{fgWex20symcond})-(\ref{fgWex21symcond}).

\end{proof}
\noindent The next lemma gives the transition pdf for a FPT of a Wiener process $W$ at
time $t_{m+1}^n$ constrained by the absorbing boundary $g^{n}$ over $[t_{m}^n, t_{m+1}^n]$ given that $W_{t_m} = x_m$.
\begin{lemma}
\label{lemmatrans}
For any $n \in \naturels_*$ and any $m \in \naturels$ such that $t_{m+1}^n\leq t_f$ we have
\begin{eqnarray}
p_{g^{n}}^W\big(t_{m+1}^n, x_{m+1} \mid t_m^n, x_m\big) & = & \Big(1 - \exp(\frac{-2(g^{n}(t_{m+1}^n)-x_{m+1})(g^{n}(t_{m}^n)-x_{m})}{\Delta_n})\Big) \nonumber \\ & & \frac{\exp (-\frac{(x_{m+1}-x_{m})^2}{\Delta_n})}{\sqrt{\pi\Delta_n}}. \label{alpha0neq}
\end{eqnarray}
\end{lemma}
\begin{proof}[Proof of Lemma \ref{lemmatrans}]
Equation (\ref{alpha0neq}) can be obtained directly from Equation (\ref{exlineareq}) in 
Lemma \ref{exlinear}.
\end{proof}
\noindent 

\begin{proof}[Proof of Lemma \ref{lemmaexpform0}]
We have that
\begin{eqnarray} \nonumber
\proba \big(T_{\alpha_{0}^n}^W \in [0,\Delta_n]\big) &= & \proba \big((T_{\alpha_{0}^n}^W > \Delta_n\big)^C)\\ \nonumber
&= & 1 - \proba \big(T_{\alpha_{0}^n}^W > \Delta_n\big)\\ \nonumber
&= & 1 - \proba \big((W_{t_1^n} - W_{t_0^n}) \in(-\infty, \alpha_{0}^n], T_{\alpha_{0}^n}^W > \Delta_n\big)\\ \nonumber
&= & 1 - \int_{-\infty}^{\alpha_{0}^n}p_{\alpha_{0}^n}^W(t_1^n, x_1 \mid t_0^n, 0) dx_1\\ \label{PTA0expression0}
&= & 1 - \int_{-\infty}^{\alpha_{0}^n} (1 - \exp(\frac{-2(\alpha_{0}^n-x_{1})\alpha_{0}^n}{\Delta_n}))\frac{\exp (-\frac{x_{1}^2}{\Delta_n})}{\sqrt{\pi\Delta_n}} dx_1
\end{eqnarray}
where we use the fact that $T_{\alpha_{0}^n}^W \geq 0$ a.s. by Definition \ref{defFPT} along with the completeness of the filtration $\mathbf{F}$ in the first equality, elementary probability facts in the second equality, the fact that $T_{\alpha_{0}^n}^W \subset \{ (W_{t_1^n} - W_{t_0^n}) \in(-\infty, \alpha_{0}^n] \}$ in the third equality, Equation (\ref{exeqpiecewiselinear2}) from Lemma \ref{expiecewiselinear} in the fourth equality and Equation (\ref{alpha0neq}) from Lemma \ref{lemmatrans} in the fifth equality. Finally, we can deduce Equation (\ref{lemmaexpform0}) by plugging Equation (\ref{alpha0n}) into Equation (\ref{PTA0expression0}).
\end{proof}
\begin{proof}[Proof of Lemma \ref{lemmaexpform1}]
We have that
\begin{eqnarray} \nonumber
\proba \big(T_{g^n}^W \in [0,\Delta_n]\big) &= & \proba \big((T_{g^n}^W > \Delta_n\big)^C)\\ \nonumber
&= & 1 - \proba \big(T_{g^n}^W > \Delta_n\big)\\ \nonumber
&= & 1 - \proba \big((W_{t_1^n} - W_{t_0^n}) \in(-\infty, g^n(t_1^n)], T_{\alpha_{0}^n}^W > \Delta_n\big)\\ \nonumber
&= & 1 - \int_{-\infty}^{g^n(t_1^n)}p_{g^n}^W(t_1^n, x_1 \mid t_0^n, 0) dx_1\\ \nonumber
&= & 1 - \int_{-\infty}^{g^n(t_1^n)} (1 - \exp(\frac{-2(g^n(t_1^n)-x_{1})g^n(t_0^n)}{\Delta_n}))\frac{\exp (-\frac{x_{1}^2}{\Delta_n})}{\sqrt{\pi\Delta_n}} dx_1\\ \label{PTA0expression}
&= & 1 - \int_{-\infty}^{\alpha_0^n + \alpha_1^n \Delta_n} (1 - \exp(\frac{-2(\alpha_0^n + \alpha_1^n \Delta_n -x_{1})\alpha_0^n)}{\Delta_n}))\frac{\exp (-\frac{x_{1}^2}{\Delta_n})}{\sqrt{\pi\Delta_n}} dx_1
\end{eqnarray}
where we use the fact that $T_{g^n}^W \geq 0$ a.s. by Definition \ref{defFPT} along with the completeness of the filtration $\mathbf{F}$ in the first equality, elementary probability facts in the second equality, the fact that $T_{g^n}^W \subset \{ (W_{t_1^n} - W_{t_0^n}) \in(-\infty, g^n(t_1^n)] \}$ in the third equality, Equation (\ref{exeqpiecewiselinear2}) from Lemma \ref{expiecewiselinear} in the fourth equality, Equation (\ref{alpha0neq}) from Lemma \ref{lemmatrans} in the fifth equality and Equations (\ref{seq1})-(\ref{seq2}) in the sixth equality. Finally, we can deduce Equation (\ref{lemmaexpform0}) by plugging Equation (\ref{alpha0n}) into Equation (\ref{PTA0expression}).
\end{proof}
\begin{proof}[Proof of Lemma \ref{lemmaa1n}]
We have that
\begin{eqnarray}
\nonumber \proba \big(T_{g^{n}}^W \in [\Delta_n,2\Delta_n]\big) &=& \proba\left((\Tau_{g^{n}}^W < \Delta_n, \Tau_{g^{n}}^W > 2 \Delta_n)^C \right)\\ \nonumber
&= & 1 - \proba\left(\Tau_{g^{n}}^W < \Delta_n, \Tau_{g^{n}}^W > 2 \Delta_n \right)\\ \nonumber
&= & 1 - \proba\left(\Tau_{g^{n}}^W < \Delta_n \right) - \proba\left(\Tau_{g^{n}}^W > 2 \Delta_n \right)\\ \nonumber
&= & 1 - \proba\left(0 \leq \Tau_{g^{n}}^W < \Delta_n \right) - \proba\left(\Tau_{g^{n}}^W > 2 \Delta_n \right)\\\label{pgexp}
&= & 1 - \int_{0}^{\Delta_n}f(s)ds - \proba\left(\Tau_{g^{n}}^W > 2 \Delta_n \right),
\end{eqnarray}
where we use elementary probability facts in the first and second equalities, the fact that $\{\Tau_{g^{n}}^W < \Delta_n\}$ and  $\{\Tau_{g^{n}}^W > 2 \Delta_n \}$ are disjoint events in the third equality, the fact that $T_{g^n}^W \geq 0$ a.s. by Definition \ref{defFPT} along with the completeness of the filtration $\mathbf{F}$ in the fourth equality, Equation (\ref{alpha0n}) in the fifth equality. Also, we have that
\begin{eqnarray} 
\nonumber \proba\left(\Tau_{g^{n}}^W > 2 \Delta_n \right)  &= & \proba \big((W_{t_1^n} - W_{t_0^n}) \in(-\infty, g^n(t_1^n)], (W_{t_2^n} - W_{t_0^n}) \in(-\infty, g^n(t_2^n)], \\ \nonumber
& &\Tau_{g^{n}}^W > 2 \Delta_n\big)\\ \nonumber
&= & \int_{-\infty}^{g^n(t_1^n)} \int_{-\infty}^{g^n(t_2^n)} p_{g^{n}}^W(t_1^n, x_1, t_2^n, x_2 \mid t_0^n, 0) dx_1 dx_2\\ \nonumber
&= & \int_{-\infty}^{g^n(t_1^n)} \int_{-\infty}^{g^n(t_2^n)}  p_{g^{n}}^W(t_2^n, x_2 \mid t_1^n, x_1) p_{g^{n}}^W(t_1^n, x_1 \mid t_0^n, 0) dx_1 dx_2\\ \nonumber
&= & \int_{-\infty}^{g^n(t_1^n)} \int_{-\infty}^{g^n(t_2^n)} (1 \\ \nonumber & & - \exp(\frac{-2(g^n(t_2^n)-x_{2})(g^n(t_1^n)-x_{1})}{\Delta_n})) \frac{\exp (-\frac{(x_{2}-x_{1})^2}{\Delta_n})}{\sqrt{\pi\Delta_n}}\\ \nonumber & &   
(1 - \exp(\frac{-2(g^{n}(t_1^n)-x_{1})g^{n}(t_0^n)}{\Delta_n})) \frac{\exp (-\frac{x_{1}^2}{\Delta_n})}{\sqrt{\pi\Delta_n}} dx_1 dx_2\\ \label{PT2exp}
&= & \int_{-\infty}^{\alpha_0^n+\alpha_1^n \Delta_n} \int_{-\infty}^{\alpha_0^n+(\alpha_1^n + \alpha_2^n) \Delta_n} (1 - \exp(\\ \nonumber & & \frac{-2(\alpha_0^n+(\alpha_1^n + \alpha_2^n) \Delta_n-x_{2})(\alpha_0^n+\alpha_1^n \Delta_n-x_{1})}{\Delta_n})) \frac{\exp (-\frac{(x_{2}-x_{1})^2}{\Delta_n})}{\sqrt{\pi\Delta_n}}\\ \nonumber & &   
(1 - \exp(\frac{-2(\alpha_0^n+\alpha_1^n \Delta_n-x_{1})\alpha_0^n}{\Delta_n})) \frac{\exp (-\frac{x_{1}^2}{\Delta_n})}{\sqrt{\pi\Delta_n}} dx_1 dx_2,
\end{eqnarray}
where we use the fact that $$\{\Tau_{g^{n}}^W > 2 \Delta_n\} \subset \{ (W_{t_1^n} - W_{t_0^n}) \in(-\infty, g^n(t_1^n)], (W_{t_2^n} - W_{t_0^n}) \in(-\infty, g^n(t_2^n)] \}$$ 
in the first equality, Equation (\ref{exeqpiecewiselinear2}) in the second equality, Equation (\ref{pgWfact}) in the third equality, Equation (\ref{alpha0neq}) from Lemma \ref{lemmatrans} in the fourth equality and Equations (\ref{seq1})-(\ref{seq2}) in the fifth equality. Finally, we can deduce Equation (\ref{alpha1neq2sim}) by plugging Equation (\ref{PT2exp}) and Equation (\ref{pintf}) into Equation (\ref{pgexp}).
\end{proof}

\begin{proof}[Proof of Lemma \ref{lemmaexpformm}]
We have that
\begin{eqnarray}
\nonumber  \proba \big(T_{g^{n}}^W \in [m\Delta_n,(m+1)\Delta_n]\big) &=& \proba\left((\Tau_{g^{n}}^W < m \Delta_n, \Tau_{g^{n}}^W > (m+1) \Delta_n)^C \right)\\ \nonumber
&= & 1 - \proba\left(\Tau_{g^{n}}^W < m\Delta_n, \Tau_{g^{n}}^W > (m+1) \Delta_n \right)\\ \nonumber
&= & 1 - \proba\left(\Tau_{g^{n}}^W < m\Delta_n \right) - \proba\left(\Tau_{g^{n}}^W > (m+1) \Delta_n \right)\\ \nonumber
&= & 1 - \proba\left(0 \leq \Tau_{g^{n}}^W < m \Delta_n \right) - \proba\left(\Tau_{g^{n}}^W > (m+1) \Delta_n \right)\\\label{pgexpgen}
&= & 1 - \int_{0}^{m\Delta_n}f(s)ds - \proba\left(\Tau_{g^{n}}^W > (m+1) \Delta_n \right),
\end{eqnarray}
where we use elementary probability facts in the first and second equalities, the fact that $\{\Tau_{g^{n}}^W < m\Delta_n\}$ and  $\{\Tau_{g^{n}}^W > (m+1) \Delta_n \}$ are disjoint events in the third equality, the fact that $T_{g^n}^W \geq 0$ a.s. by Definition \ref{defFPT} along with the completeness of the filtration $\mathbf{F}$ in the fourth equality, Equations (\ref{alpha0n})-(\ref{pintf}) in the fifth equality. Also, we have that
\begin{eqnarray} 
\nonumber \proba\left(\Tau_{g^{n}}^W > (m+1) \Delta_n \right)  &= & \proba \big((W_{t_1^n} - W_{t_0^n}) \in(-\infty, g^n(t_1^n)], (W_{t_2^n} - W_{t_0^n}) \in(-\infty, g^n(t_2^n)], \ldots, \\ \nonumber & &  (W_{t_{m+1}^n} - W_{t_0^n}) \in(-\infty, g^n(t_{m+1}^n)], \Tau_{g^{n}}^W > (m+1) \Delta_n\big)\\ \nonumber
&= & \int_{-\infty}^{g^n(t_1^n)} \int_{-\infty}^{g^n(t_2^n)} \ldots \int_{-\infty}^{g^{n}(t_{m+1}^n)} p_{g^{n}}^W(t_1^n, x_1, t_2^n, x_2 , \ldots, t_{m+1}^n, x_{m+1} \mid \\ \nonumber & &  t_0^n, 0) dx_1 dx_2 \ldots dx_{m+1}\\ \nonumber
&= & \int_{-\infty}^{g^n(t_1^n)} \int_{-\infty}^{g^n(t_2^n)} \ldots \int_{-\infty}^{g^{n}(t_{m+1}^n)} \prod_{k=0}^{m} p_{g^{n}}^W(t_{k+1}^n, x_{k+1}  \mid t_k^n, x_k)\\ \nonumber & &  dx_1 dx_2\ldots dx_{m+1}\\ \nonumber
&= & \int_{-\infty}^{g^n(t_1^n)} \int_{-\infty}^{g^n(t_2^n)} \ldots \int_{-\infty}^{g^{n}(t_{m+1}^n)} \prod_{k=0}^{m} (1 - \\ \nonumber & & \exp(\frac{-2(g^{n}(t_{k+1}^n)-x_{k+1})(g^{n}(t_k^n)-x_{k})}{\Delta_n})) \\ \nonumber  & &\frac{\exp (-\frac{(x_{k+1}-x_{k})^2}{\Delta_n})}{\sqrt{\pi\Delta_n}} dx_1 dx_2\ldots dx_{m+1},\\ \nonumber
&= & \int_{-\infty}^{\alpha_0^n+\alpha_1^n \Delta_n} \int_{-\infty}^{\alpha_0^n+(\alpha_1^n + \alpha_2^n) \Delta_n} \ldots \int_{-\infty}^{\alpha_0^n+\sum_{k=1}^{m+1}\alpha_k^n \Delta_n} \prod_{k=0}^{m} (1 - \\ \nonumber & & \exp(\frac{-2(\alpha_0^n+\sum_{i=1}^{k+1}\alpha_i^n \Delta_n-x_{k+1})(\alpha_0^n+\sum_{i=1}^{k}\alpha_i^n \Delta_n-x_{k})}{\Delta_n})) \\ \label{PT2expgen}  & &\frac{\exp (-\frac{(x_{k+1}-x_{k})^2}{\Delta_n})}{\sqrt{\pi\Delta_n}} dx_1 dx_2\ldots dx_{m+1},
\end{eqnarray}
where we use the fact that
$$\{ \Tau_{g^{n}}^W > (m+1) \Delta_n \} \subset \{ (W_{t_1^n} - W_{t_0^n}) \in(-\infty, g^n(t_1^n)], (W_{t_2^n} - W_{t_0^n}) \in(-\infty, g^n(t_2^n)], \ldots,$$ 
$$ (W_{t_{m+1}^n} - W_{t_0^n}) \in(-\infty, g^n(t_{m+1}^n)] \},$$ 
in the first equality, Equation (\ref{exeqpiecewiselinear2}) in the second equality, Equation (\ref{pgWfact}) in the third equality, Equation (\ref{alpha0neq}) from Lemma \ref{lemmatrans} in the fourth equality and Equations (\ref{seq1})-(\ref{seq2}) in the fifth equality. Finally, we can deduce Equation (\ref{alphamneq2sim}) by plugging Equation (\ref{PT2expgen}) and Equation (\ref{alpha0n}) into Equation (\ref{pgexpgen}).
\end{proof}
\noindent The next lemma will be useful in showing the existence and unicity of $\alpha_{0}^{n}$, i.e., in the proof of Proposition \ref{propexistence}. This basically states that the probability of the FPT started at time $t_0^n$ to a constant boundary on $[t_0^n,t_1^n]$, i.e., $G_{0}^n$ or $H_{0}^n$, is a strictly decreasing bijection from $\reels^+_*$ to $(0,1)$.
\begin{lemma}
\label{lemmazi} 
For any $n \in \naturels_*$ we have that $G_{0}^n$ and $H_{0}^n$ are continuous and strictly decreasing bijections from $\reels^+_*$ to $(0,1)$.
\end{lemma}
\begin{proof}
From Equation (\ref{defF0n}), we can see that $G_{0}^n$ is continuous, and admits a derivative which is negative for any $n \in \naturels_*$. Thus we have that $G_{0}^n$ is strictly decreasing. We also have that $G_{0}^n (\alpha) \rightarrow 1$ as $\alpha \rightarrow 0$ and $G_{0}^n (\alpha) \rightarrow 0$ as $\alpha \rightarrow \infty$, thus $G_{0}^n$ is continuous a bijection from $\reels^+_*$ to $(0,1)$.
We can prove the case $H_0^n$ with the same arguments.
\end{proof}
\noindent The next lemma is the counterpart of Lemma \ref{lemmazi} when considering $G_{m}^{n}$ and $H_{m}^{n}$ for any $n \in \naturels_*$ and any $m \in \naturels_*$ such that $t_{m}^n\leq t_f$. 
\begin{lemma}
\label{lemmaR}
For any $n \in \naturels_*$ and any $m \in \naturels$ such that $t_{m}^n\leq t_f$ we have that $G_m^n$ and $H_m^n$ are continuous and strictly decreasing bijections from $\reels$ to $(0,\int_{m\Delta_n}^{+\infty}f(s)ds)$.
\end{lemma}
\begin{proof}
From Equations (\ref{defF1n})-(\ref{defF2n})-(\ref{defFmn}), we can see that $G_{m}^n$ is continuous, and admits a derivative which is negative for any $n \in \naturels_*$ and any $m \in \naturels_*$ such that $t_{m}^n\leq t_f$. Thus we have that $G_{m}^n$ is strictly decreasing. We also have that $G_{m}^n (\alpha) \rightarrow \int_{m \Delta_n}^{+\infty}f(s)ds$ as $\alpha \rightarrow - \infty$ and $G_{m}^n (\alpha) \rightarrow 0$ as $\alpha \rightarrow \infty$, thus $G_{m}^n$ is continuous a bijection $\reels$ to $(0,\int_{m \Delta_n}^{+\infty}f(s)ds)$.
We can prove the case $H_m^n$ with the same arguments. 
\end{proof}
The following lemma shows the almost everywhere positivity of $f$ when we assume that \textbf{Assumption [A]} holds.
\begin{lemma}
\label{lemmaass}
We assume that \textbf{Assumption [A]} holds. Then, we have that $f$ is positive on $\reels^+$ almost everywhere.
\end{lemma}
\begin{proof}[Proof of Lemma \ref{lemmaass}]
To prove Lemma \ref{lemmaass}, it is sufficient by Borel arguments to prove that 
\begin{eqnarray}
\label{proof23}
0 < \int_{m \Delta_n}^{(m+1) \Delta_n}f(s)ds  
\end{eqnarray}
for any $n \in \naturels_*$ and by induction on $m \in \naturels$ such that $t_m^n\leq t_f$. We start with the $m=0$ case, i.e., 
\begin{eqnarray}
0 < \int_{0}^{\Delta_n}f(s)ds. \label{proof230}
\end{eqnarray}
We define the maximum of the absolute boundary $g$ on $[t_0,t_f]$ as $$g_+ := \underset{t \in [t_0,t_f]}{\sup} \mid g(t) \mid.$$ By \textbf{Assumption [A]}, we have that $g$ is continuous on $[t_0,t_f]$, and since it is a compact space it implies that $g_+ < \infty$. By Definition \ref{defFPT}, we can deduce that $\Tau_{g}^Z \leq \Tau_{g_+}^Z$ a.s.. Thus, we can deduce that
\begin{eqnarray}
\proba(\Tau_{g_+}^Z \in [0,\Delta_n]) \leq \proba ( \Tau_{g}^Z \in [0,\Delta_n]).
\end{eqnarray}
Since $\proba(\Tau_{g_+}^Z \in [0,\Delta_n]) = G_{0}^n(g_+)$ or $\proba(\Tau_{g_+}^Z \in [0,\Delta_n]) = H_{0}^n(g_+)$, we obtain by Lemma \ref{lemmazi} that $\proba(\Tau_{g_+}^Z \in [0,\Delta_n]) > 0$. Then, we can deduce Equation (\ref{proof230}) since $f$ is equal to the density of $\Tau_{g}^Z$ by Equation (\ref{PgZdef}) and Equation (\ref{ffg}). The $m>0$ case  follows since by \textbf{Assumption [A]}, we have that $g$ admits a derivative almost everywhere which implies that $\proba ( \Tau_{g}^Z \in [m \Delta_n,(m+1)\Delta_n]) >0$.

\end{proof}

\begin{proof}[Proof of Proposition \ref{propexistence}] For any $n \in \naturels_*$, we prove Proposition \ref{propexistence} by induction on $m \in \naturels$ such that $t_m^n\leq t_f$. We start with the $m=0$ case, i.e., we show that $\alpha_{0}^n \in \reels_*^+$ is well-defined. By Lemma \ref{lemmaass} along with \textbf{Assumption [A]} we can deduce that
\begin{eqnarray}
0 < \int_{0}^{\Delta_n}f(s)ds < 1. \label{proof210217}
\end{eqnarray}
From Expression (\ref{proof210217}) and Lemma \ref{lemmaexpform0}, we can then deduce that
\begin{eqnarray}
\label{proof0331c}
0 < G_0^{n} (\alpha_{0}^{n}) < 1 \text{ and } 0 < H_0^{n} (\alpha_{0}^{n}) < 1.
\end{eqnarray}
Finally, an application of the intermediate value theorem together with Lemma \ref{lemmazi} and Expression (\ref{proof0331c}) provides the existence and uniqueness of $\alpha_{0}^{n} \in \reels_*^+$. We consider now the $m>0$ case, i.e., we show that $\alpha_{m}^{n} \in \reels$ is well-defined. 
By Lemma \ref{lemmaass} along with \textbf{Assumption [A]} we get
\begin{eqnarray}
\label{proof220301}
0 < \int_{m \Delta_n}^{(m+1) \Delta_n}f(s)ds < \int_{m \Delta_n}^{+ \infty}f(s)ds. 
\end{eqnarray}
From Expression (\ref{proof220301}) and Lemmas \ref{lemmaexpform1}-\ref{lemmaa1n}-\ref{lemmaexpformm}, we can deduce that
\begin{eqnarray}
\label{proof0401c}
0 < G_m^n(\alpha_{m}^{n}) < \int_{m \Delta_n}^{+ \infty}f(s)ds \text{ and } 0 < H_m^n(\alpha_{m}^{n}) < \int_{m \Delta_n}^{+ \infty}f(s)ds.
\end{eqnarray}
To conclude, an application of the intermediate value theorem along with Lemma \ref{lemmaR} and Equation (\ref{proof0401c}) provides the existence and uniqueness of $\alpha_{m}^{n} \in \reels$.
\end{proof}

\section{Proof of approximation convergence}
\label{existence}
In this section, we show that a subsequence of the approximation uniformly converges to the boundary when the length of each interval of linear approximation goes to 0 asymptotically. The proof goes in two steps. First, we show that the approximation uniformly converges to some boundary $\tilde{g} \in \mathcal{G}$ using Arzel\`a-Ascoli theorem on any compact space $[t_0,t_f]$. Second, we show that $\tilde{g}(t) = g(t)$ for any $t \in [t_0, t_f]$. In what follows, we give the definition of uniform boundedness and uniform equicontinuity, and the Arzel\`a-Ascoli theorem.
\begin{definition}
The sequence $g^{n} \in \mathcal{G}_{PL}^{n}$ defined on the interval $[t_0, t_f]$ is uniformly bounded if there is a constant number $M > 0$ such that
\begin{eqnarray}
\label{aa1}
\sup_{t \in [t_0, t_f], n\in \naturels^*} \left|g^{n}(t)\right|\leq M.
\end{eqnarray}
\end{definition}
\begin{definition}
The sequence $g^{n} \in \mathcal{G}_{PL}^{n}$ defined on the interval $[t_0, t_f]$ is uniformly equicontinuous if, for every $\varepsilon > 0$, there exists a $\delta > 0$ such that
\begin{eqnarray}
\label{aa2}
\sup_{t,s \in [t_0, t_f], |t - s| < \delta, n\in \naturels^*} \left|\left|g^{n}(t)-g^{n}(s)\right|\right| \leq \varepsilon.
\end{eqnarray}
\end{definition}
\begin{theorem}[Arzel\`a-Ascoli theorem] 
\label{ArzelaAscolitheorem}
If the sequence $g^{n} \in \mathcal{G}_{PL}^{n}$ defined on the interval $[t_0, t_f]$ is uniformly bounded and uniformly equicontinuous, then there exists a subsequence which converges uniformly to some $\tilde{g} \in \mathcal{G}$ defined on the interval $[t_0, t_f]$.
\end{theorem}

\smallskip
In the following proposition, we show that if we assume that the $\alpha_{m}^{n}$ are uniformly bounded, then the sequence $g^n$ is uniformly bounded and uniformly equicontinuous.
\begin{proposition}
\label{arzasci}
We assume that \textbf{Assumption [A]} and \textbf{Assumption [B]} hold. If we also assume that the $\alpha_{m}^{n}$ are uniformly bounded, i.e., 
\begin{eqnarray}
\label{Kalphamn}
\sup_{\underset{m=0,\cdots,2^n}{n \in \naturels}} \left| \alpha_{m}^{n} \right| \leq K,
\end{eqnarray} 
the sequence $g^n$ is uniformly bounded and uniformly equicontinuous, i.e., it satifies Equations (\ref{aa1})-(\ref{aa2}).
\end{proposition}

\begin{proof}
We start with the proof of Equation (\ref{aa1}). By algebraic manipulation, we can rewrite Equations (\ref{seq1})-(\ref{seq2}) as
\begin{eqnarray}
\label{seq3} 
g^{n}(u) & = & \alpha_{0}^n + \Delta_n \sum_{i=1}^{m} \alpha_{i}^n + \alpha_{m+1}^n (u - t_m^n) , u \in (t_m^n,t_{m+1}^n], m \in \naturels \text{ s.t. } t_{m+1}^n\leq t_f.
\end{eqnarray}
We obtain that for $u \in (t_m^n,t_{m+1}^n], m \in \naturels \text{ such that } t_{m+1}^n\leq t_f$ that
\begin{eqnarray*}
\left|g^{n}(u) \right|& \leq & \left|\alpha_{0}^n \right| + \Delta_n \sum_{i=1}^{m} \left|\alpha_{i}^n \right|+ \left|\alpha_{m+1}^n \right|(u - t_m^n)\\
& \leq & \left|\alpha_{0}^n \right| + \Delta_n \sum_{i=1}^{m+1} \left|\alpha_{i}^n \right|\\
& \leq & \left|\alpha_{0}^n \right|+ \Delta_n \sum_{i=1}^{2^n} \left|\alpha_{i}^n \right|\\
& \leq &\left|\alpha_{0}^n \right| + (t_f - t_0) \sup_{\underset{i=1,\cdots,2^n}{n \in \naturels}} \left|\alpha_{i}^n \right|\\
& \leq & (1 + (t_f - t_0)) \sup_{\underset{i=0,\cdots,2^n}{n \in \naturels}} \left|\alpha_{i}^n \right|\\
& \leq & (1 + (t_f - t_0)) K,
\end{eqnarray*}
where we use the triangular inequality in the first inequality, the fact that $u \in (t_m^n,t_{m+1}^n]$ in the second inequality, \textbf{Assumption [B]} in the third equality, the definition of $\Delta_n$ in the fourth equality and Equation (\ref{Kalphamn}) in the last inequality. We have thus shown that Equation (\ref{Kalphamn}) $\implies$ Equation (\ref{aa1}). We now prove Equation (\ref{aa2}). We consider any arbitrarily small $\varepsilon > 0$. Accordingly, we set 
\begin{eqnarray}
\label{delta} \delta = \frac{\varepsilon}{2K}.
\end{eqnarray}
 For any $t \in [t_0,t_f]$, we define the corresponding $m_t^n$ such that 
$t \in [t_{m_t^n}^n,t_{m_t^n+1}^n]$.
From Equation (\ref{seq3}), we can deduce that
\begin{eqnarray}
\label{proof20230820}
g^{n}(t) & = & \alpha_{0}^n + \Delta_n \sum_{i=1}^{m_t^n} \alpha_{i}^n + \alpha_{m_t^n+1}^n (t - t_{m_t^n}^n ).
\end{eqnarray}
Thus, for any $t_0 \leq s \leq t \leq t_f$ such that 
\begin{eqnarray}
\label{delta0}
|t - s| < \delta.
\end{eqnarray}
We have that
\begin{eqnarray*}
\left| g^{n}(t) - g^{n}(s) \right| & = & \Big| \alpha_{0}^n + \Delta_n \sum_{i=1}^{m_t^n} \alpha_{i}^n + \alpha_{m_t^n+1}^n (t - t_{m_t^n}^n ) - (\alpha_{0}^n + \Delta_n \sum_{i=1}^{m_s^n} \alpha_{i}^n + \alpha_{m_s^n+1}^n (s - t_{m_s^n}^n )) \Big|\\
& = & \Big| \Delta_n \sum_{i=m_s^n}^{m_t^n} \alpha_{i}^n + \alpha_{m_t^n+1}^n (t - t_{m_t^n}^n) - \alpha_{m_s^n+1}^n (s - t_{m_s^n}^n) \Big|\\
& \leq & \left| t-s \right| \sup_{\underset{i=0,\cdots,2^n}{n \in \naturels}} \left|\alpha_{i}^n \right|\\
& \leq & K \left| t-s \right|,
\\
& \leq & \varepsilon,
\end{eqnarray*}
where we use Equation (\ref{proof20230820}) in the first equality, algebraic manipulation in the second equality and the first equality, Equation (\ref{Kalphamn}) in the second inequality, Equation (\ref{delta}) and Expression (\ref{delta0}) in the last inequality. We have thus shown that Equation (\ref{Kalphamn}) $\implies$ Equation (\ref{aa2}).
\end{proof}
\noindent In the following proposition, we show that if we assume that \textbf{Assumption [A]}, \textbf{Assumption [B]} and \textbf{Assumption [C]} hold, then we have that the $\alpha_{m}^{n}$ are uniformly bounded.
\begin{proposition}
\label{propassB}
We assume that \textbf{Assumption [A]}, \textbf{Assumption [B]} and \textbf{Assumption [C]} hold. Then, we have that the $\alpha_{m}^{n}$ are uniformly bounded, i.e., 
Equation (\ref{Kalphamn}) is satisfied.
\end{proposition}
\begin{proof}
We define the bound as
\begin{eqnarray}
\label{boundK}
K = \sup (\underset{t \in [t_0,t_f]}{\sup} \mid g(t) \mid, \underset{t \in [t_0,t_f]}{\sup} \mid g'(t) \mid).
\end{eqnarray}
Since we have that $t_f$ is finite by \textbf{Assumption [B]} and that $g$ is continuous on the compact space $[t_0,t_f]$ by \textbf{Assumption [A]}, we can deduce that $$\underset{t \in [t_0,t_f]}{\sup} \mid g(t) \mid < \infty.$$ We can also obtain by \textbf{Assumption [C]} that $$\underset{t \in [t_0,t_f]}{\sup} \mid g'(t) \mid < \infty.$$ Thus, we can deduce that $K < \infty$. Moreover, $K$ does not depend on $n$ or $m$ by definition. Then, to prove Proposition \ref{propassB} it is sufficient to show that Equation (\ref{Kalphamn}) is satisfied with K defined in Equation (\ref{boundK}). For any $n \in \naturels_*$, we consider a proof by induction on $m \in \naturels$ such that $t_m^n\leq t_f$. We start with the case $m=0$, i.e., we show that $\alpha_{0}^n \leq K$. By Definition \ref{defFPT}, we can deduce that $\Tau_{g}^Z \leq \Tau_{g_+}^Z$ a.s.. Thus, we can deduce that
\begin{eqnarray*}
\proba(\Tau_{g_+}^Z \in [0,\Delta_n]) \leq \proba ( \Tau_{g}^Z \in [0,\Delta_n]).
\end{eqnarray*}
By Equations (\ref{PgZdef})-(\ref{fZgt}), the above inequality can be reexpressed as 
\begin{eqnarray*}
\proba(\Tau_{g_+}^Z \in [0,\Delta_n])  \leq \int_{0}^{\Delta_{n}} f_{g}^Z (s) ds.
\end{eqnarray*}
By Equation (\ref{ffg}), the above inequality can be reexpressed as  
\begin{eqnarray*}
\proba(\Tau_{g_+}^Z \in [0,\Delta_n])  \leq \int_{0}^{\Delta_{n}} f (s) ds.
\end{eqnarray*}
By Equation (\ref{alpha0n}), the above inequality can be reexpressed as  
\begin{eqnarray*}
\proba(\Tau_{g_+}^Z \in [0,\Delta_n])  \leq \proba(\Tau_{\alpha_0^n}^Z \in [0,\Delta_n]).
\end{eqnarray*}
By Equation (\ref{defF0n}) and Lemma \ref{lemmaexpform0}, the above inequality can be reexpressed as 
\begin{eqnarray*}
G_{0}^n(g_+) \leq G_{0}^n(\alpha_0^n),
\end{eqnarray*}
or $H_{0}^n(g_+) \leq H_{0}^n(\alpha_0^n)$. Since we have that $G_{0}^n$ and $H_{0}^n$ are continuous and strictly decreasing bijections from $\reels^+_*$ to $(0,1)$  by Lemma \ref{lemmazi}, we can deduce that $\alpha_{0}^n \leq g_+$ which implies $\alpha_{0}^n \leq K$. We consider now the case $m=1$, i.e., we show that $\mid \alpha_{1}^n \mid \leq K$. We define the maximum of the absolute boundary derivative $g'$ on $[t_0,t_f]$ as $$g_+' := \underset{t \in [t_0,t_f]}{\sup} \mid g'(t) \mid.$$ For $t \geq t_0$, we define the linear boundary started at $g(t_0)$ with trend $g_+'$ and  $- g_+'$ as respectively $\overline{g}(t) = g(t_0) + g_+'(t-t_0)$ and $\underline{g}(t) = g(t_0) - g_+'(t-t_0)$. By Definition \ref{defFPT}, we can deduce that $\Tau_{\underline{g}}^Z \leq \Tau_{g}^Z \leq \Tau_{\overline{g}}^Z$ a.s.. Thus, we can deduce that
\begin{eqnarray*}
\proba(\Tau_{\overline{g}}^Z \in [0,\Delta_n]) \leq \proba ( \Tau_{g}^Z \in [0,\Delta_n]) \leq \proba(\Tau_{\underline{g}}^Z \in [0,\Delta_n]).
\end{eqnarray*}
By Equations (\ref{PgZdef})-(\ref{fZgt}), the above inequalities can be reexpressed as 
\begin{eqnarray*}
\proba(\Tau_{\overline{g}}^Z \in [0,\Delta_n])  \leq \int_{0}^{\Delta_{n}} f_{g}^Z (s) ds \leq \proba(\Tau_{\underline{g}}^Z \in [0,\Delta_n]).
\end{eqnarray*}
By Equation (\ref{ffg}), the above inequalities can be reexpressed as  
\begin{eqnarray*}
\proba(\Tau_{\overline{g}}^Z \in [0,\Delta_n]) \leq \int_{0}^{\Delta_{n}} f (s) ds\leq \proba(\Tau_{\underline{g}}^Z \in [0,\Delta_n]).
\end{eqnarray*}
By Equation (\ref{pintf}), the above inequalities can be reexpressed as  
\begin{eqnarray*}
\proba(\Tau_{\overline{g}}^Z \in [0,\Delta_n])  \leq \proba(\Tau_{g^n}^Z \in [0,\Delta_n])\leq \proba(\Tau_{\underline{g}}^Z \in [0,\Delta_n]).
\end{eqnarray*}
By Equation (\ref{defF1n}) and Lemma \ref{lemmaexpform1}, the above inequalities can be reexpressed as 
\begin{eqnarray*}
G_{1}^n(g_+') \leq G_{1}^n(\alpha_1^n) \leq G_{1}^n(-g_+'),
\end{eqnarray*}
or $H_{1}^n(g_+') \leq H_{1}^n(\alpha_1^n) \leq H_{1}^n(-g_+')$. Since we have that $G_1^n$ and $H_1^n$ are continuous and strictly decreasing bijections from $\reels$ to $(0,\int_{ \Delta_n}^{+\infty}f(s)ds)$  by Lemma \ref{lemmaR}, we can deduce that $\mid \alpha_{1}^n \mid \leq g_+'$ which implies $\mid \alpha_{1}^n \mid \leq K$. By Equation (\ref{defF1n}) and Lemma \ref{lemmaexpform1}, the above inequalities can be reexpressed as 
\begin{eqnarray*}
G_{1}^n(g_+') \leq G_{1}^n(\alpha_1^n) \leq G_{1}^n(-g_+'),
\end{eqnarray*}
or $H_{1}^n(g_+') \leq H_{1}^n(\alpha_1^n) \leq H_{1}^n(-g_+')$. Since we have that $G_1^n$ and $H_1^n$ are continuous and strictly decreasing bijections from $\reels$ to $(0,\int_{ \Delta_n}^{+\infty}f(s)ds)$  by Lemma \ref{lemmaR}, we can deduce that $\mid \alpha_{1}^n \mid \leq g_+'$ which implies $\mid \alpha_{1}^n \mid \leq K$. We consider now the case $m=2$, i.e., we show that $\mid \alpha_{2}^n \mid \leq K$. For $t \geq t_0$, we define the boundary which is equal to $g$ on $[t_0^n,t_1^n]$ and linear with trend $g_+'$ and  $- g_+'$ for $t \geq t_1^n$ as respectively 
\begin{eqnarray*}
\overline{g}(t) & = & g(t) \text{ for any } t \in [t_0^n,t_1^n]\\\overline{g}(t) & = & g(t_1^n)+ g_+'(t-t_1^n) \text{ for any } t \geq t_1^n
\end{eqnarray*}
and
\begin{eqnarray*}
\underline{g}(t) & = & g(t) \text{ for any } t \in [t_0^n,t_1^n]\\\underline{g}(t) & = & g(t_1^n) - g_+'(t-t_1^n) \text{ for any } t \geq t_1^n.
\end{eqnarray*}
By Definition \ref{defFPT}, we can deduce that $\Tau_{\underline{g}}^Z \leq \Tau_{g}^Z \leq \Tau_{\overline{g}}^Z$ a.s.. Thus, we can deduce that
\begin{eqnarray*}
\proba(\Tau_{\overline{g}}^Z \in [\Delta_n,2\Delta_n]) \leq \proba ( \Tau_{g}^Z \in [\Delta_n,2\Delta_n]) \leq \proba(\Tau_{\underline{g}}^Z \in [\Delta_n,2\Delta_n]).
\end{eqnarray*}
By Equations (\ref{PgZdef})-(\ref{fZgt}), the above inequalities can be reexpressed as 
\begin{eqnarray*}
\proba(\Tau_{\overline{g}}^Z \in [\Delta_n,2\Delta_n])  \leq \int_{\Delta_n}^{2\Delta_{n}} f_{g}^Z (s) ds \leq \proba(\Tau_{\underline{g}}^Z \in [\Delta_n,2\Delta_n]).
\end{eqnarray*}
By Equation (\ref{ffg}), the above inequalities can be reexpressed as  
\begin{eqnarray*}
\proba(\Tau_{\overline{g}}^Z \in [\Delta_n,2\Delta_n]) \leq \int_{\Delta_n}^{2\Delta_{n}} f (s) ds\leq \proba(\Tau_{\underline{g}}^Z \in [\Delta_n,2\Delta_n]).
\end{eqnarray*}
By Equation (\ref{pintf}), the above inequalities can be reexpressed as  
\begin{eqnarray*}
\proba(\Tau_{\overline{g}}^Z \in [\Delta_n,2\Delta_n])  \leq \proba(\Tau_{g^n}^Z \in [\Delta_n,2\Delta_n])\leq \proba(\Tau_{\underline{g}}^Z \in [\Delta_n,2\Delta_n]).
\end{eqnarray*}
By Equation (\ref{defF2n}) and Lemma \ref{lemmaa1n}, the above inequalities can be reexpressed as 
\begin{eqnarray*}
G_{2}^n(g_+') \leq G_{2}^n(\alpha_2^n) \leq G_{2}^n(-g_+'),
\end{eqnarray*}
or $H_{2}^n(g_+') \leq H_{2}^n(\alpha_2^n) \leq H_{2}^n(-g_+')$. Since we have that $G_2^n$ and $H_2^n$ are continuous and strictly decreasing bijections from $\reels$ to $(0,\int_{ 2\Delta_n}^{+\infty}f(s)ds)$  by Lemma \ref{lemmaR}, we can deduce that $\mid \alpha_{2}^n \mid \leq g_+'$ which implies $\mid \alpha_{2}^n \mid \leq K$. The case with $m>2$ follows with the same arguments.
\end{proof}
The following corollary is an application of Arzel\`a-Ascoli theorem.
\begin{corollary}
\label{corollaryAA}
We assume that \textbf{Assumption [A]}, \textbf{Assumption [B]} and \textbf{Assumption [C]} hold. Then, there exists a subsequence $g^{n_k}$ of $g^{n}$ which converges uniformly to some $\tilde{g} \in \mathcal{G}$ defined on the interval $[t_0, t_f]$.
\end{corollary}
\begin{proof}
This is an application of Theorem \ref{ArzelaAscolitheorem} (Arzel\`a-Ascoli theorem) along with Proposition \ref{arzasci} and Proposition \ref{propassB}.
\end{proof}

The following lemma gives a.s. convergence of $T_{h^{n}}^Z \mathbf{1}_{\{ T_{h^{n}}^Z \leq t_f - t_0 \}}$ to $T_{h}^Z \mathbf{1}_{\{ T_h^Z \leq t_f - t_0 \}}$ when $h^{n}$ converges uniformly to $h$ on $[t_0,t_f]$. For the proof of convergence, we only require the convergence in distribution.
\begin{lemma}
\label{lemma0}
For any sequence $h^n \in \mathcal{G}^n$ which converges uniformly on $[t_0,t_f]$ to some $h\in \mathcal{G}$ satisfying \textbf{Assumption A}, we have that $\Tau_{h^{n}}^Z \mathbf{1}_{\{ \Tau_{h^{n}}^Z \leq t_f - t_0 \}}$ converges a.s. to $\Tau_{h}^Z \mathbf{1}_{\{ \Tau_h^Z \leq t_f - t_0 \}}$. As a by-product, we deduce that $\Tau_{h^{n}}^Z$ converges in distribution to $\Tau_{h}^Z$ on $[t_0,t_f]$. 
\end{lemma}
\begin{proof}
To prove that $\Tau_{h^{n}}^Z \mathbf{1}_{\{ \Tau_{h^{n}}^Z \leq t_f - t_0 \}}$ converges a.s. to $\Tau_{h}^Z \mathbf{1}_{\{ \Tau_h^Z \leq t_f - t_0 \}}$, it is sufficient to show that for any arbitrarily small $\epsilon > 0$ there exists $N_\epsilon \in \naturels$ such that for any $n \in \naturels_*$ with $n \geq N_\epsilon$ we have a.s.
 \begin{eqnarray}
\label{proof0430}
 \left|\Tau_{h^{n}}^Z \mathbf{1}_{\{ \Tau_{h^{n}}^Z \leq t_f - t_0 \}} - \Tau_{h}^Z \mathbf{1}_{\{ \Tau_h^Z \leq t_f - t_0 \}} \right| \leq \epsilon.
\end{eqnarray}
As $h^{n}$ converges uniformly to $h$ on $[t_0,t_f]$, we have that for any $\epsilon_h > 0$, there exists $N_{\epsilon_h} \in \naturels$ such that for any $n \in \naturels_*$ with $n \geq N_{\epsilon_h}$ we have 
\begin{eqnarray}
\label{epsilonh0}
\underset{t \in [t_0,t_f]}{\sup} \mid h^{n}(t) - h(t) \mid \leq \epsilon_h.
\end{eqnarray}
We set the value of $\epsilon_h$ as 
\begin{eqnarray}
\label{epsilonh}
\epsilon_h = \frac{1}{2} \sup_{\Tau_{h}^Z \leq t \leq \Tau_{h}^Z + \epsilon \leq t_f - t_0} \left|Z_t - h (t)\right|.
\end{eqnarray}
First, we can see that $\epsilon_h$ defined in Equation (\ref{epsilonh}) is positive. Second, we have that a.s. $Z_t$ first hits $h^{n}$ on $[\Tau_{h}^Z - \epsilon, T_{h}^Z + \epsilon]$, i.e., we have shown that $\Tau_{h^{n}}^Z \in [\Tau_{h}^Z - \epsilon, \Tau_{h}^Z + \epsilon]$ whenever Equation (\ref{epsilonh0}) holds with $\epsilon_h$ from Equation (\ref{epsilonh}). Thus, we have shown Equation (\ref{proof0430}) with $N_\epsilon = N_{\epsilon_h}$.
\end{proof}
\noindent We consider a discretization length in the order $\frac{1}{2^n}$ so that we obtain that the time discretization is nested, i.e., for any $t_m^n$ and any $l \geq m$ there exists a time $t_k^l$ such that $t_m^n=t_k^l$. This is required to prove the following lemma which in turn will be used to prove that the limit of a subsequence obtained by Arzel\`a-Ascoli theorem satisfies Equation (\ref{ffg}).
\begin{lemma}
\label{corexistence}
We assume that \textbf{Assumption [A]} holds. For any $n \in \naturels_*$, any $l \in \naturels_*$ with $l \geq n$ and any $m \in \naturels \text{ s.t. } t_{m+1}^n\leq t_f$, the approximated boundary satisfies
\begin{eqnarray}
\label{propexeq2}
\proba \left(\Tau_{g^{l}}^Z \in [t_m^n,t_{m+1}^n]\right) = \int_{t_m^n}^{t_{m+1}^n}f(s)ds.
\end{eqnarray} 
\end{lemma}
\begin{proof}
For any $n \in \naturels_*$, any $l \in \naturels_*$ with $l \geq n$ and any $m \in \naturels \text{ s.t. } t_{m+1}^n\leq t_f$, we have
\begin{eqnarray*}
\proba \left(\Tau_{g^{l}}^Z \in [t_m^n,t_{m+1}^n]\right) & = & \sum_{i \in \naturels \text{ s.t. } t_{m}^n \leq t_{i}^l \leq t_{i+1}^l \leq t_{m+1}^n} \proba \left(\Tau_{g^{l}}^Z \in [t_i^l,t_{i+1}^l]\right)\\
& = & \sum_{i \in \naturels \text{ s.t. } t_{m}^n \leq t_{i}^l \leq t_{i+1}^l \leq t_{m+1}^n} \int_{t_{i}^l}^{t_{i+1}^l}f(s)ds\\
& = &\int_{t_{m}^n}^{t_{m+1}^n}f(s)ds,
\end{eqnarray*} 
where we use the fact that $[t_m^n,t_{m+1}^n]=\bigcup_{i \in \naturels \text{ s.t. } t_{m}^n \leq t_{i}^l \leq t_{i+1}^l \leq t_{m+1}^n}  [t_i^l,t_{i+1}^l]$ since the time discretization is nested in the first equality, and Equations (\ref{alpha0n})-(\ref{pintf}) in the second equality.
\end{proof}
We provide in what follows the proof of the main result of our paper, which shows that a subsequence of the new approximation uniformly converges to the boundary when the length of each interval of linear approximation goes to 0 asymptotically. This proof is based on application of previously obtained results and shows that $\tilde{g}(t) = g(t)$ for any $t \in [t_0, t_f]$.
\begin{proof}[Proof of Theorem \ref{main}]
By Corollary \ref{corollaryAA} along with \textbf{Assumption [A]}-\textbf{Assumption [B]}-\textbf{Assumption [C]}, there exists a subsequence $g^{n_k}$ of $g^{n}$ which converges uniformly to some $\tilde{g} \in \mathcal{G}$ defined on the interval $[t_0, t_f]$. We first show that the density of $f_{\tilde{g}}^Z(t)=f(t)$ for any $t \in [0,t_f-t_0]$. By Borel arguments, it is sufficient to show that for any $p \in \naturels_*$ and $k=0,\ldots,2^{p}-1$ we have
\begin{eqnarray}
\label{proof0}
\proba \big(\Tau_{\tilde{g}}^Z \in [k \Delta_p,(k+1) \Delta_p]\big) = \int_{k \Delta_p}^{(k+1) \Delta_p}f(s)ds.
\end{eqnarray}
We have that
\begin{eqnarray*}
\proba \big(\Tau_{\tilde{g}}^Z \in [k \Delta_p,(k+1) \Delta_p]\big) & = & \lim_{n \rightarrow \infty} \proba \big(\Tau_{g^{n_k}}^Z \in [k \Delta_p,(k+1) \Delta_p]\big) \\
& = & \int_{k\Delta_p}^{(k+1)\Delta_p}f(s)ds,
\end{eqnarray*}
where the first equality corresponds to the convergence in distribution of $\Tau_{g^{n_k}}^Z$ to $\Tau_{\tilde{g}}^Z$ by Lemma \ref{lemma0} along with \textbf{Assumption [A]}, and we use Lemma \ref{corexistence} in the second equality. Thus, we have shown Equation (\ref{proof0}), which implies that $f_{\tilde{g}}^Z(t)=f(t)$ for any $t \in [0,t_f-t_0]$. Since there is uniqueness of the inverse first-passage problem by the papers mentioned in the introduction,
we can deduce that $\tilde{g}(t) = g(t)$ for any $t \in [t_0, t_f]$.
\end{proof}




\newpage 

\begin{funding}
The author was supported in part by Japanese Society for the Promotion of Science Grants-in-Aid for Scientific Research (B) 23H00807 and Early-Career Scientists 20K13470. 
\end{funding}



\bibliographystyle{imsart-nameyear} 
\bibliography{biblio}       

\begin{thebibliography}{38}

\bibitem[\protect\citeauthoryear{Abbring}{2012}]{abbring2012mixed}
\begin{barticle}[author]
\bauthor{\bsnm{Abbring},~\bfnm{Jaap~H}\binits{J.~H.}}
(\byear{2012}).
\btitle{Mixed Hitting-Time Models}.
\bjournal{Econometrica}
\bvolume{80}
\bpages{783--819}.
\end{barticle}
\endbibitem

\bibitem[\protect\citeauthoryear{Abundo}{2006}]{abundo2006limit}
\begin{barticle}[author]
\bauthor{\bsnm{Abundo},~\bfnm{Mario}\binits{M.}}
(\byear{2006}).
\btitle{Limit at zero of the first-passage time density and the inverse problem
  for one-dimensional diffusions}.
\bjournal{Stochastic analysis and applications}
\bvolume{24}
\bpages{1119--1145}.
\end{barticle}
\endbibitem

\bibitem[\protect\citeauthoryear{Anderson}{1960}]{anderson1960modification}
\begin{barticle}[author]
\bauthor{\bsnm{Anderson},~\bfnm{Todd~W}\binits{T.~W.}}
(\byear{1960}).
\btitle{A modification of the sequential probability ratio test to reduce the
  sample size}.
\bjournal{The Annals of Mathematical Statistics}
\bvolume{31}
\bpages{165--197}.
\end{barticle}
\endbibitem

\bibitem[\protect\citeauthoryear{Anulova}{1981}]{anulova1981markov}
\begin{barticle}[author]
\bauthor{\bsnm{Anulova},~\bfnm{SV}\binits{S.}}
(\byear{1981}).
\btitle{On Markov stopping times with a given distribution for a Wiener
  process}.
\bjournal{Theory of Probability \& Its Applications}
\bvolume{25}
\bpages{362--366}.
\end{barticle}
\endbibitem

\bibitem[\protect\citeauthoryear{Beiglb{\"o}ck
  et~al.}{2018}]{beiglbock2018geometry}
\begin{barticle}[author]
\bauthor{\bsnm{Beiglb{\"o}ck},~\bfnm{Mathias}\binits{M.}},
  \bauthor{\bsnm{Eder},~\bfnm{Manu}\binits{M.}},
  \bauthor{\bsnm{Elgert},~\bfnm{Christiane}\binits{C.}} \AND
  \bauthor{\bsnm{Schmock},~\bfnm{Uwe}\binits{U.}}
(\byear{2018}).
\btitle{Geometry of distribution-constrained optimal stopping problems}.
\bjournal{Probability Theory and Related Fields}
\bvolume{172}
\bpages{71--101}.
\end{barticle}
\endbibitem

\bibitem[\protect\citeauthoryear{Butler and
  Huzurbazar}{1997}]{butler1997stochastic}
\begin{barticle}[author]
\bauthor{\bsnm{Butler},~\bfnm{Ronald~W}\binits{R.~W.}} \AND
  \bauthor{\bsnm{Huzurbazar},~\bfnm{Aparna~V}\binits{A.~V.}}
(\byear{1997}).
\btitle{Stochastic network models for survival analysis}.
\bjournal{Journal of the American Statistical Association}
\bvolume{92}
\bpages{246--257}.
\end{barticle}
\endbibitem

\bibitem[\protect\citeauthoryear{Chen, Chadam and
  Saunders}{2022}]{chen2022higher}
\begin{barticle}[author]
\bauthor{\bsnm{Chen},~\bfnm{Xinfu}\binits{X.}},
  \bauthor{\bsnm{Chadam},~\bfnm{John}\binits{J.}} \AND
  \bauthor{\bsnm{Saunders},~\bfnm{David}\binits{D.}}
(\byear{2022}).
\btitle{Higher-order regularity of the free boundary in the inverse
  first-passage problem}.
\bjournal{SIAM Journal on Mathematical Analysis}
\bvolume{54}
\bpages{4695--4720}.
\end{barticle}
\endbibitem

\bibitem[\protect\citeauthoryear{Chen et~al.}{2011}]{chen2011existence}
\begin{barticle}[author]
\bauthor{\bsnm{Chen},~\bfnm{Xinfu}\binits{X.}},
  \bauthor{\bsnm{Cheng},~\bfnm{Lan}\binits{L.}},
  \bauthor{\bsnm{Chadam},~\bfnm{John}\binits{J.}} \AND
  \bauthor{\bsnm{Saunders},~\bfnm{David}\binits{D.}}
(\byear{2011}).
\btitle{Existence and uniqueness of solutions to the inverse boundary crossing
  problem for diffusions}.
\bjournal{Annals of Applied Probability}
\bvolume{21}
\bpages{1663--1693}.
\end{barticle}
\endbibitem

\bibitem[\protect\citeauthoryear{Cheng et~al.}{2006}]{cheng2006analysis}
\begin{barticle}[author]
\bauthor{\bsnm{Cheng},~\bfnm{Lan}\binits{L.}},
  \bauthor{\bsnm{Chen},~\bfnm{Xinfu}\binits{X.}},
  \bauthor{\bsnm{Chadam},~\bfnm{John}\binits{J.}} \AND
  \bauthor{\bsnm{Saunders},~\bfnm{David}\binits{D.}}
(\byear{2006}).
\btitle{Analysis of an inverse first passage problem from risk management}.
\bjournal{SIAM Journal on Mathematical Analysis}
\bvolume{38}
\bpages{845--873}.
\end{barticle}
\endbibitem

\bibitem[\protect\citeauthoryear{Daniels}{1969}]{daniels1969minimum}
\begin{barticle}[author]
\bauthor{\bsnm{Daniels},~\bfnm{Henry~E}\binits{H.~E.}}
(\byear{1969}).
\btitle{The minimum of a stationary Markov process superimposed on a U-shaped
  trend}.
\bjournal{Journal of Applied Probability}
\bvolume{6}
\bpages{399--408}.
\end{barticle}
\endbibitem

\bibitem[\protect\citeauthoryear{Doob}{1949}]{doob1949heuristic}
\begin{barticle}[author]
\bauthor{\bsnm{Doob},~\bfnm{Joseph~L}\binits{J.~L.}}
(\byear{1949}).
\btitle{Heuristic approach to the Kolmogorov-Smirnov theorems}.
\bjournal{The Annals of Mathematical Statistics}
\bvolume{20}
\bpages{393--403}.
\end{barticle}
\endbibitem

\bibitem[\protect\citeauthoryear{Dudley and Gutmann}{1977}]{dudley1977stopping}
\begin{barticle}[author]
\bauthor{\bsnm{Dudley},~\bfnm{Richard~M}\binits{R.~M.}} \AND
  \bauthor{\bsnm{Gutmann},~\bfnm{Sam}\binits{S.}}
(\byear{1977}).
\btitle{Stopping times with given laws}.
\bjournal{S{\'e}minaire de Probabilit{\'e}s de Strasbourg}
\bvolume{11}
\bpages{51--58}.
\end{barticle}
\endbibitem

\bibitem[\protect\citeauthoryear{Durbin}{1971}]{durbin1971boundary}
\begin{barticle}[author]
\bauthor{\bsnm{Durbin},~\bfnm{James}\binits{J.}}
(\byear{1971}).
\btitle{Boundary-crossing probabilities for the Brownian motion and Poisson
  processes and techniques for computing the power of the Kolmogorov-Smirnov
  test}.
\bjournal{Journal of Applied Probability}
\bvolume{8}
\bpages{431--453}.
\end{barticle}
\endbibitem

\bibitem[\protect\citeauthoryear{Durbin}{1985}]{durbin1985first}
\begin{barticle}[author]
\bauthor{\bsnm{Durbin},~\bfnm{James}\binits{J.}}
(\byear{1985}).
\btitle{The first-passage density of a continuous Gaussian process to a general
  boundary}.
\bjournal{Journal of Applied Probability}
\bvolume{22}
\bpages{99--122}.
\end{barticle}
\endbibitem

\bibitem[\protect\citeauthoryear{Ekstr{\"o}m and
  Janson}{2016}]{ekstrom2016inverse}
\begin{barticle}[author]
\bauthor{\bsnm{Ekstr{\"o}m},~\bfnm{Erik}\binits{E.}} \AND
  \bauthor{\bsnm{Janson},~\bfnm{Svante}\binits{S.}}
(\byear{2016}).
\btitle{The inverse first-passage problem and optimal stopping}.
\bjournal{Annals of Applied Probability}
\bvolume{26}
\bpages{3154--3177}.
\end{barticle}
\endbibitem

\bibitem[\protect\citeauthoryear{Gut}{1974}]{gut1974moments}
\begin{barticle}[author]
\bauthor{\bsnm{Gut},~\bfnm{Allan}\binits{A.}}
(\byear{1974}).
\btitle{On the moments and limit distributions of some first passage times}.
\bjournal{The Annals of Probability}
\bvolume{2}
\bpages{277--308}.
\end{barticle}
\endbibitem

\bibitem[\protect\citeauthoryear{Jacod and Shiryaev}{2003}]{JacodLimit2003}
\begin{bbook}[author]
\bauthor{\bsnm{Jacod},~\bfnm{Jean}\binits{J.}} \AND
  \bauthor{\bsnm{Shiryaev},~\bfnm{Albert}\binits{A.}}
(\byear{2003}).
\btitle{Limit theorems for stochastic processes},
\bedition{2nd} ed.
\bpublisher{Berlin: Springer-Verlag}.
\end{bbook}
\endbibitem

\bibitem[\protect\citeauthoryear{Klump and Kolb}{2023}]{klump2023uniqueness}
\begin{barticle}[author]
\bauthor{\bsnm{Klump},~\bfnm{Alexander}\binits{A.}} \AND
  \bauthor{\bsnm{Kolb},~\bfnm{Martin}\binits{M.}}
(\byear{2023}).
\btitle{Uniqueness of the Inverse First-Passage Time Problem and the Shape of
  the Shiryaev Boundary}.
\bjournal{Theory of Probability \& Its Applications}
\bvolume{67}
\bpages{570--592}.
\end{barticle}
\endbibitem

\bibitem[\protect\citeauthoryear{Klump and Savov}{2023}]{klump2023conditions}
\begin{barticle}[author]
\bauthor{\bsnm{Klump},~\bfnm{Alexander}\binits{A.}} \AND
  \bauthor{\bsnm{Savov},~\bfnm{Mladen}\binits{M.}}
(\byear{2023}).
\btitle{Conditions for existence and uniqueness of the inverse first-passage
  time problem applicable for Levy processes and diffusions}.
\bjournal{arXiv preprint arXiv:2305.10967}.
\end{barticle}
\endbibitem

\bibitem[\protect\citeauthoryear{Lai}{1974}]{lai1974control}
\begin{barticle}[author]
\bauthor{\bsnm{Lai},~\bfnm{Tze~Leung}\binits{T.~L.}}
(\byear{1974}).
\btitle{Control charts based on weighted sums}.
\bjournal{The Annals of Statistics}
\bvolume{2}
\bpages{134--147}.
\end{barticle}
\endbibitem

\bibitem[\protect\citeauthoryear{Lai and Siegmund}{1977}]{lai1977nonlinear}
\begin{barticle}[author]
\bauthor{\bsnm{Lai},~\bfnm{T.~L.}\binits{T.~L.}} \AND
  \bauthor{\bsnm{Siegmund},~\bfnm{D.}\binits{D.}}
(\byear{1977}).
\btitle{A nonlinear renewal theory with applications to sequential analysis I}.
\bjournal{The Annals of Statistics}
\bvolume{5}
\bpages{946--954}.
\end{barticle}
\endbibitem

\bibitem[\protect\citeauthoryear{Lai and Siegmund}{1979}]{lai1979nonlinear}
\begin{barticle}[author]
\bauthor{\bsnm{Lai},~\bfnm{T.~L.}\binits{T.~L.}} \AND
  \bauthor{\bsnm{Siegmund},~\bfnm{D.}\binits{D.}}
(\byear{1979}).
\btitle{A nonlinear renewal theory with applications to sequential analysis
  II}.
\bjournal{The Annals of Statistics}
\bvolume{7}
\bpages{60--76}.
\end{barticle}
\endbibitem

\bibitem[\protect\citeauthoryear{Malmquist}{1954}]{malmquist1954certain}
\begin{barticle}[author]
\bauthor{\bsnm{Malmquist},~\bfnm{Sten}\binits{S.}}
(\byear{1954}).
\btitle{On certain confidence contours for distribution functions}.
\bjournal{The Annals of Mathematical Statistics}
\bvolume{25}
\bpages{523--533}.
\end{barticle}
\endbibitem

\bibitem[\protect\citeauthoryear{Matthews, Farewell and
  Pyke}{1985}]{matthews1985asymptotic}
\begin{barticle}[author]
\bauthor{\bsnm{Matthews},~\bfnm{DE}\binits{D.}},
  \bauthor{\bsnm{Farewell},~\bfnm{VT}\binits{V.}} \AND
  \bauthor{\bsnm{Pyke},~\bfnm{R}\binits{R.}}
(\byear{1985}).
\btitle{Asymptotic score-statistic processes and tests for constant hazard
  against a change-point alternative}.
\bjournal{The Annals of Statistics}
\bvolume{13}
\bpages{583--591}.
\end{barticle}
\endbibitem

\bibitem[\protect\citeauthoryear{Mehr and McFadden}{1965}]{mehr1965certain}
\begin{barticle}[author]
\bauthor{\bsnm{Mehr},~\bfnm{CYRUS~BOZORG}\binits{C.~B.}} \AND
  \bauthor{\bsnm{McFadden},~\bfnm{JA1998850234}\binits{J.}}
(\byear{1965}).
\btitle{Certain properties of Gaussian processes and their first-passage
  times}.
\bjournal{Journal of the Royal Statistical Society Series B: Statistical
  Methodology}
\bvolume{27}
\bpages{505--522}.
\end{barticle}
\endbibitem

\bibitem[\protect\citeauthoryear{Novikov}{1971}]{novikov1971stopping}
\begin{barticle}[author]
\bauthor{\bsnm{Novikov},~\bfnm{Aleksandr~Aleksandrovich}\binits{A.~A.}}
(\byear{1971}).
\btitle{On stopping times for a Wiener process}.
\bjournal{Theory of Probability \& Its Applications}
\bvolume{16}
\bpages{449--456}.
\end{barticle}
\endbibitem

\bibitem[\protect\citeauthoryear{Novikov, Frishling and
  Kordzakhia}{1999}]{novikov1999approximations}
\begin{barticle}[author]
\bauthor{\bsnm{Novikov},~\bfnm{Alex}\binits{A.}},
  \bauthor{\bsnm{Frishling},~\bfnm{Volf}\binits{V.}} \AND
  \bauthor{\bsnm{Kordzakhia},~\bfnm{Nino}\binits{N.}}
(\byear{1999}).
\btitle{Approximations of boundary crossing probabilities for a Brownian
  motion}.
\bjournal{Journal of Applied Probability}
\bvolume{39}
\bpages{1019--1030}.
\end{barticle}
\endbibitem

\bibitem[\protect\citeauthoryear{Peskir}{2002a}]{peskir2002integral}
\begin{barticle}[author]
\bauthor{\bsnm{Peskir},~\bfnm{Goran}\binits{G.}}
(\byear{2002}a).
\btitle{On integral equations arising in the first-passage problem for Brownian
  motion}.
\bjournal{The Journal of Integral Equations and Applications}
\bvolume{14}
\bpages{397--423}.
\end{barticle}
\endbibitem

\bibitem[\protect\citeauthoryear{Peskir}{2002b}]{peskir2002limit}
\begin{barticle}[author]
\bauthor{\bsnm{Peskir},~\bfnm{Goran}\binits{G.}}
(\byear{2002}b).
\btitle{Limit at zero of the Brownian first-passage density}.
\bjournal{Probability Theory and Related Fields}
\bvolume{124}
\bpages{100--111}.
\end{barticle}
\endbibitem

\bibitem[\protect\citeauthoryear{Renault, Van~der Heijden and
  Werker}{2014}]{renault2014dynamic}
\begin{barticle}[author]
\bauthor{\bsnm{Renault},~\bfnm{Eric}\binits{E.}}, \bauthor{\bparticle{Van~der}
  \bsnm{Heijden},~\bfnm{Thijs}\binits{T.}} \AND
  \bauthor{\bsnm{Werker},~\bfnm{Bas~JM}\binits{B.~J.}}
(\byear{2014}).
\btitle{The dynamic mixed hitting-time model for multiple transaction prices
  and times}.
\bjournal{Journal of Econometrics}
\bvolume{180}
\bpages{233--250}.
\end{barticle}
\endbibitem

\bibitem[\protect\citeauthoryear{Salminen}{1988}]{salminen1988first}
\begin{barticle}[author]
\bauthor{\bsnm{Salminen},~\bfnm{Paavo}\binits{P.}}
(\byear{1988}).
\btitle{On the first hitting time and the last exit time for a Brownian motion
  to/from a moving boundary}.
\bjournal{Advances in Applied Probability}
\bvolume{20}
\bpages{411--426}.
\end{barticle}
\endbibitem

\bibitem[\protect\citeauthoryear{Siegmund}{1986}]{siegmund1986boundary}
\begin{barticle}[author]
\bauthor{\bsnm{Siegmund},~\bfnm{David}\binits{D.}}
(\byear{1986}).
\btitle{Boundary crossing probabilities and statistical applications}.
\bjournal{The Annals of Statistics}
\bvolume{14}
\bpages{361--404}.
\end{barticle}
\endbibitem

\bibitem[\protect\citeauthoryear{Song and Zipkin}{2011}]{song2011approximation}
\begin{barticle}[author]
\bauthor{\bsnm{Song},~\bfnm{Jing-Sheng}\binits{J.-S.}} \AND
  \bauthor{\bsnm{Zipkin},~\bfnm{Paul}\binits{P.}}
(\byear{2011}).
\btitle{An approximation for the inverse first passage time problem}.
\bjournal{Advances in Applied Probability}
\bvolume{43}
\bpages{264--275}.
\end{barticle}
\endbibitem

\bibitem[\protect\citeauthoryear{Strassen}{1967}]{strassen1967almost}
\begin{binproceedings}[author]
\bauthor{\bsnm{Strassen},~\bfnm{Volker}\binits{V.}}
(\byear{1967}).
\btitle{Almost sure behavior of sums of independent random variables and
  martingales}.
In \bbooktitle{Proc. Fifth Berkeley Sympos. Math. Statist. and Probability
  (Berkeley, Calif., 1965/66)}
\bvolume{2}
\bpages{315--343}.
\end{binproceedings}
\endbibitem

\bibitem[\protect\citeauthoryear{Wang and
  P{\"o}tzelberger}{1997}]{wang1997boundary}
\begin{barticle}[author]
\bauthor{\bsnm{Wang},~\bfnm{Liqun}\binits{L.}} \AND
  \bauthor{\bsnm{P{\"o}tzelberger},~\bfnm{Klaus}\binits{K.}}
(\byear{1997}).
\btitle{Boundary crossing probability for Brownian motion and general
  boundaries}.
\bjournal{Journal of Applied Probability}
\bvolume{34}
\bpages{54--65}.
\end{barticle}
\endbibitem

\bibitem[\protect\citeauthoryear{Woodroofe}{1976}]{woodroofe1976renewal}
\begin{barticle}[author]
\bauthor{\bsnm{Woodroofe},~\bfnm{Michael}\binits{M.}}
(\byear{1976}).
\btitle{A renewal theorem for curved boundaries and moments of first passage
  times}.
\bjournal{The Annals of Probability}
\bvolume{4}
\bpages{67--80}.
\end{barticle}
\endbibitem

\bibitem[\protect\citeauthoryear{Woodroofe}{1977}]{woodroofe1977second}
\begin{barticle}[author]
\bauthor{\bsnm{Woodroofe},~\bfnm{Michael}\binits{M.}}
(\byear{1977}).
\btitle{Second order approximations for sequential point and interval
  estimation}.
\bjournal{The Annals of Statistics}
\bvolume{7}
\bpages{984--995}.
\end{barticle}
\endbibitem

\bibitem[\protect\citeauthoryear{Zucca and Sacerdote}{2009}]{zucca2009inverse}
\begin{barticle}[author]
\bauthor{\bsnm{Zucca},~\bfnm{Cristina}\binits{C.}} \AND
  \bauthor{\bsnm{Sacerdote},~\bfnm{Laura}\binits{L.}}
(\byear{2009}).
\btitle{On the inverse first-passage-time problem for a Wiener process}.
\bjournal{The Annals of Applied Probability}
\bvolume{19}
\bpages{1319--1346}.
\end{barticle}
\endbibitem

\end{thebibliography}

\end{document}